\documentclass[a4paper]{easychair}

\title{G\"{o}del-McKinsey-Tarski and (not quite) Blok-Esakia for Heyting-Lewis Implication}

\date{}

\author{
Jim de Groot\inst{1}
\and
    Tadeusz Litak\inst{2}
\and
   Dirk Pattinson\inst{1}
}

\institute{%
  School of Computing
  The Australian National University\\
  \email{jim.degroot@anu.edu.au, dirk.pattinson@anu.edu.au}
\and%
  Chair For Theoretical Computer Science
  FAU Erlangen-Nuremberg \\
   \email{tadeusz.litak@fau.de}
 }
 
\authorrunning{}
\titlerunning{}

\usepackage[backend=bibtex,bibstyle=alphabetic,citestyle=alphabetic,firstinits=true,maxnames=6,
    maxalphanames=3]{biblatex}

\usepackage{hyperref}

\usepackage{amssymb}
\usepackage{amsmath}
\usepackage[shortlabels]{enumitem}
\usepackage{stmaryrd}
\usepackage{tikz-cd}
  \usetikzlibrary{arrows,arrows.meta}
  \tikzcdset{arrow style=tikz, diagrams={>=stealth'}}
\usepackage[mathcal]{euscript}
\usepackage[scr=dutchcal]{mathalfa}
\usepackage{fancyvrb}

\newcommand{\todo}[1]{\textcolor{red}{\textbf{\texttt{#1}}}}

\definecolor{uuxgreen}{cmyk}{1,0,0.75,0}
\definecolor{uuxred}{cmyk}{0.2,1,0.9,0.1}
\definecolor{uuyblue}  {cmyk}{0.9,0.55,0,0}

\makeatletter
\@ifpackageloaded{fixme}{}{\usepackage[final]{fixme}} 
\makeatother
\FXRegisterAuthor{av}{aav}{AV}
\FXRegisterAuthor{tmil}{atmil}{TL minor}
\FXRegisterAuthor{tl}{atl}{TL}
\FXRegisterAuthor{jim}{ajim}{JG} 
\FXRegisterAuthor{dirk}{adirk}{DP} 
\definecolor{mydarkgreen}{rgb}{0,0.34,0}
\newcommand{\tlnt}[1]{\tmilnote[inline,marginclue]{\textcolor{mydarkgreen}{#1}}}

\newcommand{\jim}[1]{\jimnote[inline,marginclue]{\textcolor{uuyblue}{\texttt{#1}}}}
\newcommand{\tadeusz}[1]{\tlnote[inline,marginclue]{\textcolor{purple}{\texttt{#1}}}}
\newcommand{\jg}[1]{\jim{#1}}

\newcommand{\placeholder}[1]{\textcolor{orange}{\textbf{\textsc{placeholder:}} #1}}

\newcommand{\rfse}[1]{\S~\ref{sec:#1}}
\newcommand{\rfsse}[1]{\S~\ref{subsec:#1}}

\usepackage{amsthm}
  \theoremstyle{definition}
  \swapnumbers
    \newtheorem{para}{}[section]
    \newtheorem{definition}[para]{Definition}
    \newtheorem{example}[para]{Example}
    \newtheorem{examples}[para]{Examples}
    \newtheorem{remark}[para]{Remark}
  \theoremstyle{theorem}
    \newtheorem{lemma}[para]{Lemma}
    \newtheorem{corollary}[para]{Corollary}
    \newtheorem{theorem}[para]{Theorem}
    \newtheorem{proposition}[para]{Proposition}
    
    \newtheorem*{question}{Question}

  \DeclareFontFamily{U}{mathc}{}
  \DeclareFontShape{U}{mathc}{m}{it}%
  {<->s*[1.03] mathc10}{}
  \DeclareMathAlphabet{\mathfun}{U}{mathc}{m}{it}
  
\newcommand{\mc}[1]{\mathcal{#1}}

\newcommand{\mf}[1]{\mathfrak{#1}}
\renewcommand{\sf}[1]{\mathsf{#1}}
\newcommand{\ms}[1]{\mathscr{#1}}
\newcommand{\cat}[1]{\mathtt{#1}}   
\newcommand{\fun}[1]{\mathfun{#1}}  
\newcommand{\lan}[1]{\mathcal{#1}}  
\newcommand{\alg}[1]{\mathscr{#1}}  

\newcommand{\lna}[1]{\ensuremath{\mathsf{#1}}} 
    
\renewcommand{\phi}{\varphi}
\renewcommand{\epsilon}{\varepsilon}
\renewcommand{\tilde}[1]{\widetilde{#1}}
\renewcommand{\hat}[1]{\widehat{#1}}
\newcommand{\ov}[1]{\overline{#1}}

\DeclareSymbolFont{symbolsC}{U}{txsyc}{m}{n}
\DeclareMathSymbol{\sto}{\mathrel}{symbolsC}{74}

\newcommand{\jff}{\quad\text{iff}\quad}

\newcommand{\tto}{\sto}
\mathchardef\hyphen="2D
\DeclareMathOperator{\op}{op}
\DeclareMathOperator{\Prop}{At}
\newcommand{\sFrm}{{\sto\!\hyphen\cat{Frm}}}
\newcommand{\GFrm}{\cat{G\hyphen Frm}}
\newcommand{\DFrm}{\cat{D\hyphen Frm}}
\newcommand{\sMod}{{\sto\!\hyphen\cat{Mod}}}
\newcommand{\iAMod}{\cat{iA\hyphen Mod}}
\DeclareMathOperator{\Subf}{Subf}
\newcommand{\undto}{\mathrel{\mkern1mu\underline{\mkern-1mu \to\mkern-2mu}\mkern2mu }}
\newcommand{\undsto}{\mathrel{\mkern1mu\underline{\mkern-1mu \sto\mkern-1mu}\mkern1mu }}

\newcommand{\hsigma}{\hat{\sigma}}
\newcommand{\hrho}{\hat{\rho}}

\newcommand{\llb}{\llbracket}
\newcommand{\rrb}{\rrbracket}

\newcommand{\eqc}[1]{\lfloor #1 \rfloor}

\newcommand{\bro}{\textup{(}}
\newcommand{\brc}{\textup{)}}

\DeclareMathSymbol{\vlbox}{\mathbin}{symbolsC}{107}
\DeclareMathSymbol{\slashbox}{\mathbin}{symbolsC}{108}

\newcommand{\ibx}{\mathbin{[\mathsf{i}]}}
\newcommand{\mbx}{\mathbin{[\mathsf{m}]}}
\newcommand{\idm}{\mathbin{\langle\mathsf{i}\rangle}}
\newcommand{\mdm}{\mathbin{\langle\mathsf{m}\rangle}}

\newcommand{\Boxi}{\Box_{\sf{i}}}
\newcommand{\Boxm}{\Box_{\sf{m}}}
\newcommand{\Diamondi}{\Diamond_{\sf{i}}}

\newcommand{\Ri}{R_{\sf{i}}}
\newcommand{\Rm}{R_{\sf{m}}}

\newcommand{\Grzi}{\hyperlink{ax:Grz}{\lna{Grz_i}}}

\newcommand{\HLAs}{\cat{HLAs}}

\newcommand{\To}{\Rightarrow}

\newcommand{\iA}{\lna{iA}}
\newcommand{\HL}{\lna{HL}}

\newcommand{\iAm}{\lna{iA^-}}

\newcommand{\sflb}{\lna{S4BHL}}
\newcommand{\sfk}{\lna{S4K}}
\newcommand{\lanBM}{\lan{L}_{\sf{i}, \sf{m}}}
\newcommand{\know}{\Box} 

\newcommand{\axiom}[1]{%
\renewcommand{\labelenumi}{$\mathsf{#1}$}
\renewcommand{\labelenumii}{$\mathsf{#1}$}
  \renewcommand{\theenumi}{$\mathsf{#1}$}%
  \renewcommand{\theenumii}{$\mathsf{#1}$}%
  \item
}
\newcommand{\axref}[1]{\mathbin{\text{\ref{#1}}}}  

\newcommand{\sBox}{{\scriptscriptstyle \Box}}

\newcommand{\ipr}{\lna{i}\text{-}}
\newcommand{\loga}[1]{\ipr\ref{ax:#1}}
\newcommand{\logabox}{\lna{\ipr Box}\,}
\newcommand{\logm}[1]{\ipr\ref{ax:#1}$^-$}

\newcommand{\iKb}{\lna{IntK_\sBox}}

\newcommand{\axst}[1]{\ensuremath{\mathsf{#1}}}

\newcommand{\wzmax}{\sigma^{\axst{mix}}}

\newcommand{\dist}{\kern 0.9pt}

\DeclareMathOperator{\possible}{\text{\tikz[scale=.5ex/1cm,baseline=-.6ex,rotate=45,line width=.1ex]{
                            \draw (-1,-1) rectangle (1,1);}}\dist}
\DeclareMathOperator{\necessary}{\text{\tikz[scale=.6ex/1cm,baseline=-.6ex,line width=.1ex]{
                            \draw (-1,-1) rectangle (1,1);}}\dist}

\renewcommand{\Box}{\necessary}
\renewcommand{\Diamond}{\possible}                                                        

\begin{document}

\maketitle

\begin{abstract}
Heyting-Lewis Logic is the extension of intuitionistic propositional
logic with a strict implication connective that satisfies the constructive counterparts of axioms for strict implication provable in classical modal logics. 
Variants of this logic are surprisingly widespread: they 
appear as Curry-Howard correspondents of (simple type theory
extended with) Haskell-style arrows, 
in preservativity logic of
Heyting arithmetic, in the proof theory of guarded (co)recursion, and in
the generalization of intuitionistic epistemic logic.

Heyting-Lewis Logic can be interpreted in 
intuitionistic Kripke frames extended with a binary relation to account
for strict implication. 
We use this semantics to define descriptive
frames (generalisations of Esakia spaces),
and establish a categorical duality between the algebraic
interpretation and the
frame semantics. 
We then adapt a transformation by Wolter and Zakharyaschev to
translate Heyting-Lewis Logic to classical modal logic
with two unary operators. 
This allows us to classical results to obtain 
the finite model property and decidability for a large family of
Heyting-Lewis logics.
%
\end{abstract}



%

\section{Introduction}


\noindent
Modern modal logic was invented by C.~I.~Lewis
\cite{Lewis14,Lewis18,Lewis32:book} as the theory of \emph{strict
implication} $\tto$. Lewis assumed a classical
propositional base and definability of $\tto$ in terms of unary
modal operators\footnote{Curiously, Lewis was not using $\Box$ as a
primitive, so in fact his intuitionistically problematic definition
of $\phi \tto \psi$ was $\neg\Diamond(\phi \wedge \neg\psi)$. See
\cite[App. D]{LitVis18} for an account of problems caused by Lewis'
use of a Boolean propositional base, namely trivialization
\cite{Lewis20,lewisincap} of his original system
\cite{Lewis14,Lewis18}, which in turn finally lead him to propose
systems \lna{S1}--\lna{S3} \cite[App. 2]{Lewis32:book} as successive
``lines of retreat''  \cite{Parry70}. Lewis considered \lna{S4} and
\lna{S5}, suggested by Becker \cite{Becker30}, too strong to provide
a proper account of strict implication \cite[p. 502]{Lewis32:book}
and appeared frustrated with later development of modal logic
\cite[\S\ 2.1]{LitVis18}. Yet, despite his supportive attitude
towards non-classical logics, he seems to have mentioned Brouwer
only once (favourably) \cite{Lewis32:monist}, and does not appear to
have ever referred to, or even be familiar with subsequent work of
Kolmogorov, Heyting or Glivenko \cite[\S\ 2.2]{LitVis18}.} as
$\Box(\phi \to \psi)$. Consequently, while implication-like
connectives are intensively studied in other areas (relevance,
substructural, counterfactual and conditional
logics), some of which in fact are of modal origin \cite[Ch.~6]{Mares04},
modal logic in a narrow sense gradually came to focus
mostly on unary boxes and diamonds, with Lewis' original strict implication $\tto$ falling
into disuse.

  Recently, Litak and Visser \cite{LitVis18} investigated $\tto$
  over an intuitionistic rather than classical propositional base,
  using intuitionistic Kripke frames with an additional \emph{binary} relation
  interpret strict implication.
  While a $\Box$-modality can be obtained from $\sto$
  via $\Box\phi := \top \sto \phi$, strict implication is \emph{not}
  definable from $\Box$.
  The constructive implication $\tto$ 
   was first studied in the context of
  \emph{preservativity} for theories over Heyting Arithmetic
  \lna{HA} (\rfse{HA}). Arrows in functional programming
  \cite{Hughes00:scp} yield another important class of Heyting-Lewis
  implications via the Curry-Howard correspondence (\rfse{haskell}).
  Yet another nontrivial variant of $\tto$ 
  arises when one generalizes Artemov and Protopopescu's \cite{ArtPro16}  approach to
  intuitionistic epistemic logic 
  (\rfse{iele}).  Even where
  $\tto$ \emph{is} reducible to $\Box$ in terms of theoremhood, it
  can still be a more useful primitive.  This
  has been argued in the proof theory of guarded (co)recursion
  (Exm.~\ref{ex:sl} in \rfse{degenerate}; see also 
  \cite[\S\ 7.2]{LitVis18}). 


While Kripke semantics has obvious benefits, it does not provide a
fully global completeness theorem for \emph{arbitrary} extensions of
the minimal Heyting-Lewis system $\iA$ (\rfse{axioms}), even in the limitative modal or superintuitionistic cases.
In contrast, the systems presented here are
amenable to algebraic semantics, called \emph{Heyting-Lewis algebras}. These are obtained by fusing Heyting algebras with so-called weak Heyting algebras \cite{CelJan05} over the shared lattice reduct  (\rfse{algebras}). 
To combine advantages of algebraic and relational semantics, one
typically works with dual representations of algebras
called descriptive frames.  These can sometimes be viewed as
topological spaces, analogous to Esakia spaces for 
intuitionistic propositional logic.  Like the algebraic semantics,
they give  completeness, but are often 
easier to manipulate and transform, e.g., to
prove the finite model property and decidability.
We provide suitable Heyting-Lewis dualities in
\rfse{sem}.

Intuitionistic logics with natural Kripke semantics can often be
viewed as fragments of classical modal logics determined by the same
Kripke structures. In particular, formulae of \lna{IPC} can be
identified via the G\"odel-McKinsey-Tarski translation with those
formulae of modal logic \lna{S4} where every subformula is prefixed
with $\Box$. 
A \emph{modal companion} of an intermediate logic $\Theta$ is then defined
as an \lna{S4}-logic containing the G\"odel-McKinsey-Tarski of an intuitionistic
formula $\phi$ if and only if $\phi$ is in $\Theta$.
These are useful because properties such as decidability and
Kripke completeness can be transferred from a modal companion back to
the corresponding intermediate logic.
%

  Wolter and Zakharyaschev \cite{WolZak97,WolZak98} extended the G\"odel-McKinsey-Tarski translation to translate intuitionistic unimodal logics into classical \emph{bimodal} logics. This proved  a fruitful approach, 
 enabling the use of well-developed classical metatheory in proofs of completeness, canonicity, the finite model property, and decidability results. In \rfse{embedding}, we generalise their result to the Heyting-Lewis setting.
  
  
  In order to put this transfer apparatus to good use, in \rfse{bimodal} we set out to prove the finite model property
  and decidability for classical bimodal logics. In particular, we prove
  this for (cofinal) \emph{transitive subframe logics}.
  These results can then be transferred to a large class of Heyting-Lewis
  logics (with an additional axiom enforcing that the relation interpreting
  $\tto$ is transitive).
  
  
\paragraph{Earlier version.}
%
  This is an updated version of a conference paper originally called
  ``G\"odel-McKinsey-Tarski and Blok-Esakia for Heyting-Lewis Implication''~\cite{GroLitPat21lics}.

  After publication of this paper and its corresponding
  technical report on arXiv~\cite{GroLitPat21arxiv},
  we were made aware of a mistake in one
  of the lemmas leading up to the Blok-Esakia theorem.
  Specifically, in Appendix A.4 of the technical report~\cite{GroLitPat21arxiv},
  the seventh equality in the chain of equalities
  at the bottom of page 36 does not hold.
  As a consequence, Lemma 4.18 fails, and with it so does Theorem 4.20,
  the analogue of the Blok-Esakia theorem.
  
  Our main application of the Blok-Esakia theorem was to transfer properties
  such as decidability, the finite model property and Kripke completeness
  from classical bimodal logics to Heyting-Lewis logics.
  To this end, we proved a general theorem deriving the finite model property
  for certain extensions of the classical bimodal logic
  $\lna{S4} \otimes \lna{K4}$.
  As a small mercy, this transfer from classical to intuitionistic can also
  be proven without the high-tech machinery of the Blok-Esakia theorem.
  
  In this updated version of the paper, we have removed the (incorrect)
  Blok-Esakia theorem and the lemmas leading up to it.
  Instead, we prove in Theorem~\ref{thm:reflect} that decidability,
  Kripke completeness and the finite
  model property can be transferred from any modal companion of a
  Heyting-Lewis logic $\Theta$ back to $\Theta$.
  We leave the open question whether or not a Blok-Esakia theorem for
  Heyting-Lewis logic can be proven as an interesting direction
  for future research.

\paragraph{Acknowledgement.}
  The authors are grateful to Cheng Liao for his interest in and careful
  reading of the predecessor of the current paper,
  and pointing out the mistake mentioned above.

\section{Syntax, Axioms and Examples}

  \noindent
  Define the language $\lan{L}_{\sto}$ by the grammar
  $$
    \phi ::=  p
         \mid \top
         \mid \bot
         \mid \phi \wedge \phi
         \mid \phi \vee \phi
         \mid \phi \to \phi
         \mid \phi \sto \phi,
  $$
  where $p$ ranges over some fixed set of propositional atoms $\Prop$.
  As usual, $\neg\phi := \phi \to \bot$. Furthermore, $\Box\phi := \top \tto \phi$.
  The unary connectives $\neg$ and $\Box$ bind strongest,
  next comes $\sto$, then $\wedge$ and $\vee$, and lastly $\to$.

\subsection{Axioms and Rules for Arrows} \label{sec:axioms}
  
  \noindent
  We define Heyting-Lewis Logic (the system $\iA$ following
  \cite{LitVis18})\footnote{Litak and Visser \cite{LitVis18} use the name ``Lewis arrow'' for $\tto$, which leads to names such as \lna{iA}, or to the use of \lna{a} as a subscript.}  as the extension of the intuitionistic propositional
  calculus ($\lna{IPC})$ with the axioms
   \begin{enumerate}
    \axiom{K_a} \label{ax:Ka}
                $((\phi \sto \psi) \wedge (\phi \sto \chi)) \to (\phi \sto (\psi \wedge \chi))$
    \axiom{Di}  \label{ax:Di}
                $((\phi \sto \chi) \wedge (\psi \sto \chi)) \to ((\phi \vee \psi) \sto \chi)$
    \axiom{Tr}  \label{ax:Tr}
                $((\phi \sto \psi) \wedge (\psi \sto \chi)) \to (\phi \sto \chi)$
\end{enumerate}
and 
the \emph{arrow necessitation} rule:
 \begin{enumerate}
    \axiom{N_a} \label{ax:Na}
                $\dfrac{\phi \to \psi}{\phi \sto \psi}$.                
  \end{enumerate}
 We also call $\iA$ the (base) \emph{Heyting-Lewis logic}. The system obtained by removing \ref{ax:Di} from the above axiomatization will be denoted as $\iAm$. In several important applications below, one needs to distinguish between the $\iAm$- and $\iA$-variants.
  Only the latter can be given a sound \emph{and complete} Kripke-style semantics
  (see \rfsse{relational}), and the results established in this paper generally
  require the presence of \ref{ax:Di} (cf. \rfse{conclusions}).
 
A ($\iA$-)\emph{logic} is a set of $\lan{L}_{\sto}$-formulae
containing all of the above axioms, and closed under \ref{ax:Na} and
uniform substitiution. Given $\Lambda, \Gamma \subseteq \lan{L}_{\sto}$, $\Lambda \oplus \Gamma$ denotes the smallest logic containing $\Lambda \cup \Gamma$. We write $\Lambda \oplus \{\phi\}$ as $\Lambda \oplus \phi$ and in the special case of  $\Lambda = \iA$ ($\Lambda = \iA^-$), we write $\ipr\phi$ ($\ipr\phi^-$).

\subsection{Intuitionistic Normal Modal Logics (with Box)} \label{sec:degenerate}

\noindent
One easily shows \cite{Iemhoff03:mlq,IemhoffJZ05:igpl,LitVis18} that
 the defined box is \emph{normal}: the axiom \lna{K_\sBox} and the rule \lna{N_\sBox} obtained by substituting $\top$ for $\phi$ in \ref{ax:Ka} and \ref{ax:Na}, respectively,  
are derivable in $\iAm$, just like
\[
\Box(\phi \to \psi) \to \phi \tto \psi.
\]
Thus, postulating as an axiom the opposite implication
  \begin{enumerate}[leftmargin=2.8cm]
    \axiom{Box}\label{ax:Box}
    $(\phi \sto \psi) \to \Box(\phi \to \psi)$
  \end{enumerate}
not only makes $\tto$ interdefinable with $\Box$ (and makes
\ref{ax:Di} derivable even over $\iAm$ \cite[Lem 4.4c]{LitVis18}),
but reduces the study of Heyting-Lewis logics extending \logabox 
to the study of normal modal extensions of \iKb, i.e., extensions of
$\lna{IPC}$ closed under \lna{K_\sBox}, \lna{N_\sBox} and uniform
substitutions in the language $\lan{L}_\sBox$ (replacing $\phi \tto
\psi$ with $\Box(\phi \to \psi)$). 

Conversely, this means that the Heyting-Lewis Logic of strict
implication subsumes a large
class of intuitionistic modal logics. 
In this rather broad area
\cite{Curry52:jsl,Ono77,BozDos84,Simpson94:phd,WolZak99,dePaivaGM04:jlc,Litak14:trends,Kavvos20}\footnote{Sotirov
\cite{Sotirov84:ml} claims accurately that ``there is an outburst of
interest in this topic each decade" since Fitch's 1948 paper
\cite{Fitch48}.} one sometimes includes a separate
$\Diamond$ connective undefinable in \logabox, but there are
often good reasons to work in the setting of a single $\Box$ added
to the \lna{IPC} signature, particularly when studying the
Curry-Howard-Lambek correspondence for a specific functor/type
operator. Examples of interest include: 
\begin{example}
The (monoidal) comonadic box of constructive \lna{S4} \cite{alechina-et-al-cs4,BiermanP00:sl}, obtained by extending \iKb\ with \hypertarget{ax:Fb}{}
\begin{enumerate}
  \axiom{T} \label{ax:T}
  $\phi \to \Box\phi \qquad\qquad\qquad \lna{4_{\sBox}} \;\; \Box\phi \to \Box\Box\phi$ \hfill
\end{enumerate}
is used to control staged computation \cite{DaviesP01:jacm,NanevskiP05:jfp}.
\end{example}
\begin{example}
The \emph{strength} axiom\footnote{\emph{Strength} of the functor interpreting $\Box$ in a categorical semantics of modal proofs  \cite{alechina-et-al-cs4,BiermanP00:sl,dePaivaR11,LitakPR17,Rog20} corresponds to the validity of $(\phi \wedge \Box\psi) \to \Box(\phi \wedge \psi)$, but this is derivable from \ref{ax:Sb} when $\Box$ is normal \cite[Sec. 6]{LitakPR17}. Classically, \loga{Sb} collapses to a fairly non-interesting system \cite[Rem. 25]{LitakPR17}.}:
\begin{enumerate}
    \axiom{S_\sBox} \label{ax:Sb}
   $\phi \to \Box\phi$.
\end{enumerate}
 yields $\iKb \oplus \lna{S_\sBox}$\footnote{We are abusing the $\oplus$ notation for logics in $\lan{L}_\sBox$, replacing closure under axioms and rules of $\iA$ with closure under axioms and rules of $\iKb$.}, the (inhabitation) logic of Haskell's \emph{applicative functors} (\emph{idioms}) \cite{McbrideP08:jfp}, as noted in recent references \cite{LitakPR17,Rog20} (cf. \rfse{haskell}). It has also been proposed as a minimal system of intuitionistic epistemic logic (cf. \rfse{iele}). 
 \end{example}
\begin{example}
Extending $\iKb \oplus \lna{S_\sBox}$ with 
\begin{enumerate}
\axiom{C4} \label{ax:CF}
$\Box\Box\phi \to \Box\phi$
\end{enumerate}
yields the Propositional Lax Logic \lna{PLL} \cite{Curry52:jsl,FaiMen97,Gol10}. This is known as the Curry-Howard correspondent of \emph{strong monads} \cite{Mog91,kobayashi-monads,benton-et-al-comp-types}, but has numerous other application in hardware verification \cite{FaiMen97}, access control \cite{GarPfe06}, epistemic logic \cite{ArtPro16}, or topos logic \cite{Gol10}.
\end{example}

\begin{example} \label{ex:sl}
  The strong L\"ob axiom 
  \begin{enumerate}[leftmargin=0.8cm]
    \axiom{SL_\sBox} \label{ax:SLb}
    $(\Box\phi \to \phi) \to \phi$ 
  \end{enumerate}
  entails \ref{ax:Sb} \cite[Lem. 3.2]{MilLit17}%
    \footnote{This is a ``deboxed'' version of well-known derivation of
    transitivity from the standard L\"ob axiom \cite[Thm.~18]{Boo93}.
    The above reference provides a categorical
    translation of this derivation.}
    and hence  is equivalent to
  \begin{enumerate} [leftmargin=0.78cm]
    \axiom{SL'_\sBox} \label{ax:SLbalt}
    $(\Box\phi \to \phi) \to \Box\phi$.
  \end{enumerate}
  Furthermore, the system $\iKb \oplus \axref{ax:SLb}$ is easily seen to be equivalent
  to the one obtained by extending $\iKb \oplus \axref{ax:Sb}$ with the (ordinary) L\"ob axiom
  \begin{enumerate}
    \axiom{L_\sBox} \label{ax:Lb}
    $\Box(\Box\phi \to \phi) \to \Box\phi$ 
  \end{enumerate}
  \ref{ax:SLb} has been studied in the context of extensions of Heyting
  Arithmetic with the \emph{completeness principle} \cite{Vis82} (see \rfse{HA}).
  However, computer scientists may 
  recognize it as an axiom for the \emph{modality for guarded \bro co\brc recursion} 
  \cite{Nak00,Nak01}, also known as the \emph{later operator} 
  \cite{BentonT09:tldi,BirMSS12,JaberTS12:lics,JunEA18},
  \emph{next clock tick} \cite{KrishnaswamiB11:lics,KrishnaswamiB11:icfp}
  or \emph{guardedness type constructor} \cite{AtkeyMB13:icfp}.
  Proof systems developed
  in this context often treat (definable) $\tto$ as a primitive connective 
  \cite{Nak00,AbeVez14,CloGor15}. Thus effectively they are proof systems
  for $\mathbin{\text{\loga{SLb}}} \oplus \axref{ax:Box}$ rather than for $\iKb \oplus \axref{ax:SLb}$;
  see \cite[\S\ 7.2]{LitVis18} for a detailed discussion.
\end{example}


\subsection{Haskell Arrows (with Choice)} \label{sec:haskell}

  \noindent
  Over $\iAm$, the strength axiom \ref{ax:Sb} is  equivalent to
  \begin{enumerate}
    \axiom{S_a}\label{ax:Sa}
    $(\phi \to \psi) \to (\phi \sto \psi)$ \cite[Lem 4.10]{LitVis18}.
  \end{enumerate}
  From a type-theoretic perspective, the axioms of \logm{Sa} correspond to
  inhabitation laws of Hughes arrows \cite{Hughes00:scp}, where an
  \emph{arrow} is a binary type constructor that represents
  computations, here given in Haskell notation:
  \begin{Verbatim}
    class Arrow a where
      arr   :: (b -> c) -> a b c
      (>>>) :: a b c ->  a c d -> a b d 
      first :: a b c ->  a (b, d) (c, d)
  \end{Verbatim}
\noindent
  Reading strict implication  $b \sto c$ as the type of arrows with
  domain $b$ and codomain $c$, the first function, \texttt{arr}
  stipulates that every function of type \texttt{b -> c} can be
  interpreted as a computation from \texttt{b} to \texttt{c}, which
  is precisely \ref{ax:Sa}.  The second function allows us to compose
  arrows, which is modally captured by \ref{ax:Tr}. Finally, the
  direct modal transliteration of \texttt{first} and \ref{ax:Ka}
  are inter-derivable over \lna{IPC} \cite[Lem. 4.1]{LitVis18}.
  This leaves the \ref{ax:Di} axiom that corresponds to a
  frequently used extension of arrows, the so-called \emph{arrows
  with choice} \cite{hughes-prog-with-arr}. This amounts to stipulating an extra operation 
  \begin{Verbatim}
  class Arrow arr => ArrowChoice arr where
    (|||) :: arr a c -> arr b c -> 
             arr (Either a b) c
   \end{Verbatim}
%
 One can set up the same
  correspondences that are usually exhibited between modal logic
  and type theory in a (strong) monad setting, e.g.,~give a realisability
  interpretation \cite{kobayashi-monads} of Heyting-Lewis proofs as
  functions in a type theory with a chosen notion of arrow, a
  propositions-as-types interpretation \cite{benton-et-al-comp-types}, or a
  categorical semantics \cite{alechina-et-al-cs4}.
  Related work has indeed been done for arrows \cite{Atk08,JacHH09,LinWY10,LinWY11},
  but inasmuch as we are aware, so far avoiding explicit mention of logic
  (or \ref{ax:Di}/Choice).%
   \footnote{To the best of our knowledge, the only explicit (if brief) discussion of the Curry-Howard connection between Haskell arrows and \logm{Sa} is found in Litak and Visser \cite[\S\ 7.1]{LitVis18}.} 
  
\begin{example}[Arrow-collapsing Choice]
Clearly, the trivial
  example of arrows, i.e. function spaces, are arrows with choice. 
  The above-mentioned applicative functors are another limiting case: they correspond to \emph{arrows with} \texttt{delay} \cite[Def. 5.1]{LinWY11}, an operation    
  which simply makes (the type corresponding to) \ref{ax:Box} 
   inhabited  \cite{LinWY11}. Finally, monads are equivalent to \emph{higher-order arrows}  \cite[\S\ 6]{LinWY11} or \emph{arrows with apply} \cite[\S 5.2]{hughes-prog-with-arr}, where the \texttt{apply} operation inhabits a type corresponding to one of Lewis' original axioms \cite{Lewis32:book}\cite[Rem. 7.3]{LitVis18}:
\begin{enumerate}[leftmargin=1cm]
\axiom{App_a}\label{ax:App}
    $(\phi \wedge (\phi \tto \psi)) \tto \psi$.
\end{enumerate}
The logic $\lna{PLAA} := \mathbin{\text{\loga{Sa}}} \oplus \axref{ax:App}$ allows for a decomposition of $\phi \tto \psi$ as $\phi \to \Box\psi$ \cite{LinWY11} \cite[Lem 4.17f]{LitVis18}. This also entails derivability of \lna{Di}  \cite[Lem 4.17g]{LitVis18}, which in the Haskell context was already noted in Hughes' original paper \cite[\S 5.2]{hughes-prog-with-arr}.\footnote{
 The logical perspective 
 seems to cast a light on the controversy whether arrows are ``stronger'' than applicative functors \cite{McbrideP08:jfp,LinWY11}. 
 Putting aside the general question of whether one takes as a measure of strength the capability to inhabit more types or rather to allow more distinctions, in the presence  of $\tto$, the situation is not as clear-cut as in the unary case, where (monadic) \lna{PLL} simply extends  (applicative) $\iKb \oplus \axref{ax:Sb}$. Defining arrows over the latter set of axioms via \texttt{delay} yields $\mathbin{\text{\logm{Sa}}} \oplus \axref{ax:Box}$, whereas inhabiting \texttt{apply} with \ref{ax:App} yields \lna{PLAA}. These are two incomparable systems.}
\end{example}

\begin{example}[Nontrivial Choice] 
  An example that does not trvialise $\tto$ is provided by Kleisli arrows:
  given a monad $M$, we define
  the type of Kleisli arrows over types $a$ and $b$ as arrows in the
  Kleisli category given by $M$, that is,
  $\mathsf{A}\, a\, b\, = a \to M\, b$. Just like with function spaces, this allows
  us to define $f ||| g = [f, g]$ as the co-pair. For co-Kleisli
  arrows, i.e., defining $\mathsf{A}\, a\, b = M\, a \to b$, we need to
  additionally require that the monad $M$ comes equipped with a
  distributive law over coproducts, viz.~$M(a + b) \to Ma + Mb$.
  Finally, list processors are presented as arrows with choice in 
  \cite{hughes-prog-with-arr}, where for the choice operation $|||$,
  the interleaving pattern in the output is modelled on the
  interleaving of the input.
\end{example}

\begin{example}[Arrows without Choice]
  An example of arrows that do not come equipped with choice are
  automata that transform elements of type $a$ to elements of type
  $b$ that satisfy the isomorphism $\mathsf{A}\, a\, b \cong a \to b
  \times (\mathsf{A}\, a\, b)$. Another non-example are functions on
  infinite streams: given two functions $\mathsf{Stream}\, a \to
  \mathsf{Stream} \, b$ and $\mathsf{Stream}\, b \to
  \mathsf{Stream}\, c$,
  there is no generic way to construct a function
  $\mathsf{Stream}\, (a + b) \to \mathsf{Stream}\, c$.
\end{example}


%

\subsection{Intuitionistic Epistemic Logic of Entailment} \label{sec:iele}

  \noindent
  From an (intuitionistic) epistemic logical point of view, Heyting-Lewis
  logic can be used to reason about entailment,
  interpreting ``$\phi \sto \psi$'' as 
  ``the agent knows that $\phi$ entails $\psi$''.
  This allows us to not only reason about the knowledge of an agent, but also
  about their deductive abilities.
  We recover statements about agent's knowledge via
  $\know\psi = \top \sto \psi$. 

  This idea leads to a generalisation of Artemov and Protopopescu's
  intuitionistic epistemic logic $\lna{IEL}$ \cite{ArtPro16}. We briefly discuss generalizations of two basic principles they postulate: \emph{coreflection} and \emph{intuitionistic reflection}. 
  
  An intuitionistic implication holds only if there exists a proof for it.
  As a consequence of this proof, our Heyting-Lewis agent \cite[\S\ 2.1]{ArtPro16} knows the implication: 
  \begin{center}
    Intuitionistic implication $\Rightarrow$ knowledge of implication
  \end{center}
  Syntactically, this simply means validity of the strength axiom \ref{ax:Sa} (or equivalently \ref{ax:Sb}), i.e., 
   \emph{coreflection}  $\phi \to \know\phi$ \cite{ArtPro16}.
  
  Conversely, known implications cannot be false.
  Therefore one cannot intuitionistically falsify any implication that is
  known. This gives rise to the following generalisation of \emph{intuitionistic reflection} ($\know\phi \to \neg\neg\phi$) 
 \begin{enumerate}
 \axiom{IR}  \label{ax:IR} $(\phi \sto \psi) \to \neg\neg(\phi \to \psi)$.
 \end{enumerate}
    
  One could say that knowledge of the entailment of $\psi$
  from $\phi$ prevents one from proving that $\neg\psi$ given $\phi$, i.e.,
  \begin{enumerate}[leftmargin=0.8cm]
    \axiom{IR'} \label{ax:IR2}
               $(\phi \sto \psi) \to (\phi \to \neg\neg\psi)$
  \end{enumerate}
  Since $\lna{IPC} \vdash \neg\neg(\phi \to \psi) \leftrightarrow (\phi \to \neg\neg\psi)$, both axioms are equivalent. 
%

\begin{definition}
  The \emph{intuitionistic epistemic logic of entailment}
  is given by
  $
    \lna{IELE} = \iA \oplus \axref{ax:Sa} \oplus \axref{ax:IR}.
  $
\end{definition}
\noindent
  There is no natural way to \ref{ax:Box}-collapse $\lna{IELE}$. 
  Knowledge of an entailment $\phi \sto \psi$ does not imply the existence
 of an intuitionistic proof, or knowledge thereof. 
  
  Besides, as a consequence of the strength axiom we have
  $(\phi \sto \psi) \to (\phi \to \Box\psi)$ (see \cite[Lem. 4.10b]{LitVis18}).
  Its converse, however, need not be true.
While $\phi$ may imply \emph{knowledge} of $\psi$, there is no
  reason it should entail \emph{intuitionistic truth} of $\psi$.
  Thus we do not wish to have the Hughes law collapsing arrows in
  the monadic setting (\rfse{haskell}), i.e., we do not stipulate
\begin{enumerate}[leftmargin=1cm]
\axiom{Hug}   \label{ax:Hug} $(\phi \to \Box \psi) \to (\phi \sto \psi)$.
\end{enumerate}
%
  We will see in Prop.~\ref{prop:IELE-form2} below that
  \ref{ax:Box}, \ref{ax:Hug} and the converse of \ref{ax:IR}
  are not derivable from $\lna{IELE}$.  
  Nevertheless, knowledge of entailment and intuitionistic implication are not entirely unrelated.
  Generalising \cite[Thm.~3.5(3)]{ArtPro16}, we find:

\begin{proposition}
  In $\lna{IELE}$ we have
  $
    \neg(\phi \sto \psi) \leftrightarrow \neg(\phi \to \psi).
  $
\end{proposition}
\begin{proof}
  The direction from left to right follows from \eqref{ax:Sa}.
  It follows from \eqref{ax:IR}
  that
  $$
    \neg\neg\neg(\phi \to \psi) \to \neg(\phi \sto \psi)
  $$
  and since $\neg(\phi \to \psi) \to \neg\neg\neg(\phi \to \psi)$ this
  proves the converse. 
\end{proof}

 Investigation of further meaningful epistemic axioms to add to \lna{IELE} is future research. See \rfse{conclusions} for a discussion.

\subsection{Arithmetical Interpretations} \label{sec:HA}

\noindent
There is more than one possible interpretation of $\tto$ in theories over Heyting Arithmetic (\lna{HA}) or other base systems such as  Intuitionistic Elementary Arithmetic (\lna{IEA}). 
A framework of \emph{schematic logics} accounts for this variety of interpretations \cite[\S\ 5]{LitVis18}, \cite[\S\ 4]{LitVis19}, broadly generalizing the well-known provability interpretation of the classical L\"ob logic \cite{Boo93}.
The most important arithmetical interpretation is that of $\Sigma_1$-preservativity, historically the first context in which constructive $\tto$ appeared \cite{Vis82,Vis94,Vis02,Iemhoff03:mlq,IemhoffJZ05:igpl}. This interpretation makes inclusion of not only \ref{ax:Sb}, but even \ref{ax:Di} in the base system problematic. 

 
 More broadly, one can define \emph{$\Delta$-preservativity}  for a theory $T$, where $\Delta$ is an elementary class of arithmetical 
 sentences containing $\top$. 
 First, define $A \tto_{\Delta,T} B$ as the following relation: for all $S \in \Delta$, if
$T \vdash S \to A$, then $T \vdash S \to B$. 
For each $\Delta$ and (a fixed axiomatization of) $T$, this yields a binary arithmetical predicate encoding the corresponding relation on G\"odel numbers in the language of \lna{IEA}.\footnote{One could go further and give a 4-argument predicate in the language of second-order arithmetic, parametric in both $\Delta$ and the axiomatization of $T$.} Now given a mapping $f$ from $\Prop$ to arithmetical sentences,  extend it inductively to the whole $\lan{L}_\tto$, using the G\"odel encoding and the predicate in question to interpret $\tto$. 
It is immediate to see that such an interpretation makes $\iAm$ valid for every $T$, $\Delta$ and $f$. 
 In order to ensure the validity of \ref{ax:Di} in the logic of $\Sigma_1$-preservativity of a given $T$, one needs additional conditions such as $T$'s $T$-provable closure under q-realizability, which does hold for \lna{HA} or Markov Arithmetic \lna{MA} \cite[\S\ 5.4.1]{LitVis18}, but not in all of their extensions  \cite[Ex~5.6]{LitVis19}. 

Unary $\Box$ encodes arithmetical provability under this interpretation, but while the usual L\"ob axiom \ref{ax:Lb} from \rfse{degenerate} remains in the provability logic of \lna{HA}, 
 the said logic contains many principles failing in the provability logic of Peano Arithmetic (\lna{PA}) \cite[\S\ 5.3]{LitVis18}. In fact, unlike the classical case \cite{Sol76}, a complete axiomatization of the intuitionistic $\lan{L}_\sBox$-provability logics remains elusive \cite{Vis94,Iem98,Iem01,Vis08,ArdMoj14} and one of the main motivations for studying the $\lan{L}_\tto$-logic of $\Sigma_1$-preservativity has been that this seemingly more challenging task may yield more natural axiomatic principles.

Furthermore, in the setting of constructive $\tto$, just the addition of \ref{ax:Lb} is not sufficient to derive the famous Explicit Fixpoint Theorem of provability logic \cite{Ber75,Sam76,BooSam91,Boo93,Vis05,vBe06} and one can consider several mutually incomparable  axioms which restore it \cite{LitVis19}. One of them is
 \begin{enumerate}
    \axiom{P} \label{ax:P} $(\phi \sto \psi) \to \Box(\phi \sto \psi)$
  \end{enumerate}
which ensures exactly the transitivity of the modal relation interpreting $\tto$ \cite[Th 10.1]{LitVis19}
and plays a central r\^{o}le in the finite model property and decidability results in \rfse{bimodal} below.

Finally, inasmuch as the strength axiom \ref{ax:Sb}  (or \ref{ax:Sa}) is concerned, it is obviously not valid in preservativity or provability logic of either \lna{HA} or \lna{PA}. Nevertheless, one can consider the system \lna{HA^*} which incorporates \ref{ax:Sb} as the so-called \emph{completeness principle} \cite{Vis82,dJV96,ArdMoj18} \cite[\S\ 5.4.4]{LitVis18}.

\section{Semantics and Duality} \label{sec:sem}

  \noindent
  The logics introduced above correspond to varieties of Heyting algebras
  with \emph{binary} operators.
  In particular, the algebraic semantics of $\iA$ is given by
  \emph{Heyting-Lewis algebras}, defined in \rfse{algebras} below.
  After defining these, we recall the relational
  semantics of $\iA$ (\rfsse{relational}),
  equip these with collections of
  admissible subsets to obtain general and descriptive frames (\rfsse{general})
  and prove a categorical duality between the descriptive frames 
  and Heyting-Lewis algebras (\rfsse{duality}).

\subsection{Heyting-Lewis Algebras} \label{sec:algebras}

\noindent
  The algebraic semantics of $\iA$ looks as follows:

\begin{definition}\label{def:iA-alg}
  A \emph{Heyting-Lewis algebra} or \emph{\HL-algebra} is a tuple of the
  form 
  $$
    \alg{A} := (A, \top, \bot, \wedge, \vee, \to, \tto),
  $$
  where $(A, \top, \bot, \wedge, \vee, \to)$ is a Heyting algebra
  and $(A, \top, \bot, \wedge, \vee, \tto)$ is a so-called
  \emph{weak Heyting algebra} \cite{CelJan05}, i.e.,
  ${\sto} : A \times A \to A$ is a binary operator that satisfies:
  \begin{enumerate}[label=$\mathsf{C}$\arabic*,leftmargin=1cm]
    \item \label{it:C1}
          $(a \sto b) \wedge (a \sto c) = a \sto (b \wedge c)$
    \item \label{it:C2}
          $(a \sto c) \wedge (b \sto c) = (a \vee b) \sto c$
    \item \label{it:C3}
          $(a \sto b) \wedge (b \sto c) \leq a \sto c$
    \item \label{it:C4}
          $a \sto a = \top$
  \end{enumerate}
  When no confusion is likely we will write $(A, \sto)$ and understand $A$
  to be (the set underlying) a Heyting algebra.

  An \emph{\HL-algebra morphism} from $(A, \sto)$ to $(A', \sto')$ is a Heyting
  homomorphism $h : A \to A'$ that additionally satisfies
  $h(a \sto b) = h(a) \sto' h(b)$ for all $a, b \in A$.
  We write $\HLAs$ for the category of \HL-algebras
  and \HL-algebra morphisms.
\end{definition}

  The collection of $\lan{L}_{\sto}$-formulae
  modulo provable equivalence yields an \HL-algebra (the Linden\-baum-Tarski algebra),  cf. \cite[\S\ 3.1]{LitVis19}. In fact, 
  completeness of $\iA$ and its extensions can be shown by standard
  techniques for \emph{algebraizable logics} \cite{BlokP89:ams,FontJP03a:sl,Rasiowa74:aatnl,Font06:sl}. 
  This is not only true for theoremhood, but also for theories induced 
  by the \emph{global consequence relation}.  As this is routine, we skip the details.

%

\subsection{Strict Implication Frames and Models}\label{subsec:relational}

  \noindent
  We recall the relational semantics for $\lan{L}_{\sto}$ \cite[\S~3.4.2]{Iem01},
  \cite[Definition 3.3]{LitVis18}.
  These are intuitionistic Kripke frames (i.e., posets) with an additional
  binary relation that is used to interpret the strict implication $\sto$.

\begin{definition}\label{def:sto-frm}
  A \emph{strict implication frame}, or \emph{$\sto$-frame} for short,
  is a tuple $(X, \preceq, \sqsubset)$ consisting of a poset $(X, \preceq)$ and
  a binary relation $\sqsubset$ on $X$ that satisfies for all $x, y, z \in X$,
  \begin{equation}\label{eq:p-sto}\tag{$\sto$-p}
    \text{if}\quad x \preceq y \sqsubset z \quad\text{then}\quad x \sqsubset z.
  \end{equation}
  A \emph{$\sto$-morphism} from $(X, \preceq, \sqsubset)$
  to $(X', \preceq', \sqsubset')$ is a function $f : X \to X'$
  that is bounded with respect to both relations. 
  That is, for $R \in \{ {\preceq}, {\sqsubset} \}$ and all $x, y \in X$ and $z' \in X'$: 
  \begin{enumerate}[label=$\mathsf{P}$\arabic*,leftmargin=1cm]
    \item If $x R y$ then $f(x) R' f(y)$;
    \item If $f(x) R' z'$ then $\exists z \in X$
          s.t.~$x R z$ and $f(z) = z'$;
  \end{enumerate}
  We write $\sFrm$ for the category of $\sto$-frames and -morphisms.
\end{definition}
  
  For a poset $(X, \preceq)$, let $\fun{up}(X, \preceq) = \{ a \subseteq X \mid \text{ if } x \in a \text{ and } x \preceq y \text{ then } y \in a \}$ be the collection
  of upsets of $X$.
  Recall that $\fun{up}(X, \preceq)$ can be given a Heyting algebra structure
  where top, bottom, meet and join are given by $X$, $\emptyset$, intersection
  and union, respectively. Implication is defined by
  $$
    a \undto b
      = \{ x \in X \mid \text{if } x \preceq y
                      \text{ and } y \in a \text{ then } y \in b \}.
  $$
  Likewise, $\sto$-frames give rise to \HL-algebras:
  
\begin{lemma}\label{lem:fun-plus-obj}
  Every $\sto$-frame $\mf{X} = (X, \preceq, \sqsubset)$ gives rise
  to an \HL-algebra $\mf{X}^+ = (\fun{up}(X, \preceq), \undsto)$, where
  $\undsto$ is defined by
  $$
    a \undsto b
      = \{ x \in X \mid \text{if } x \sqsubset y
                      \text{ and } y \in a \text{ then } y \in b \}.
  $$
\end{lemma}
\begin{proof}
  If $a$ and $b$ are upsets, then by \eqref{eq:p-sto} so is $a \undsto b$,
  so $\undsto$ is well defined.
  To prove that $(\fun{up}(X, \preceq), \undsto)$ is an \HL-algebra
  we need to show that it satisfies \ref{it:C1} to \ref{it:C4} from
  Definition \ref{def:iA-alg}. Each of these follows from a straightforward verification.
\end{proof}

  The algebra $\mf{X}^+$ is known as the \emph{complex algebra} of $\mf{X}$.
  
\begin{proposition}\label{prop:fun-plus}
  The assignment $(\cdot)^+$ extends to a contravariant functor
  $(\cdot)^+ : \sFrm \to \HLAs$ by setting $f^+ = f^{-1}$ for every
  $\sto$-frame morphism $f$.
\end{proposition}
\begin{proof}
  We have already seen that $(\cdot)^+$ is well defined on objects.
  If $f : (X, \preceq, \sqsubset) \to (X', \preceq', \sqsubset')$ is a
  $\sto$-frame morphism, then in particular it is a bounded morphism
  from $(X, \preceq)$ to $(X', \preceq')$ and hence
  $f^{-1} : \fun{up}(X', \preceq') \to \fun{up}(X, \preceq)$
  is a Heyting homomorphism.
  Boundedness of $f$ with respect to $\sqsubset$ entails that
  $f^{-1}(a' \undsto' b') = f^{-1}(a') \undsto f^{-1}(b)$.
  Functoriality is straightforward.
\end{proof}

  Thus we can choose upsets of $(X, \preceq)$ as the intepretants
  of a $\sto$-frame $(X, \preceq, \sqsubset)$ and define a
  $\sto$-model as follows.

\begin{definition}
  A \emph{valuation} for a $\sto$-frame $(X, \preceq, \sqsubset)$ is a
  function $V : \Prop \to \fun{up}(X, \preceq)$ that assigns to each propositional
  variable $p$ an upset of $(X, \preceq)$.
  A \emph{$\sto$-model} is a tuple $(\mf{X}, V)$ of a $\sto$-frame
  $\mf{X} = (X, \preceq, \sqsubset)$ and a valuation $V$ for $\mf{X}$.
  
  The set of states of a model $\mf{M} = (\mf{X}, V)$ satisfying an
  $\lan{L}_{\sto}$-formula $\phi$ is called the \emph{truth set} of $\phi$.
  It is denoted by $\llb \phi \rrb^{\mf{M}}$ and defined recursively by
  $\llb p \rrb^{\mf{M}} = V(p)$,
  $\llb \top \rrb^{\mf{M}} = X$, $\llb \bot \rrb^{\mf{M}} = \emptyset$, and
  \begin{align*}
    \llb \phi \wedge \psi \rrb^{\mf{M}}
      &= \llb \phi \rrb^{\mf{M}} \cap \llb \psi \rrb^{\mf{M}}
      &\; \llb \phi \to \psi \rrb^{\mf{M}}
      &= \llb \phi \rrb^{\mf{M}} \undto \llb \psi \rrb^{\mf{M}} \\
    \llb \phi \vee \psi \rrb^{\mf{M}}
      &= \llb \phi \rrb^{\mf{M}} \cup \llb \psi \rrb^{\mf{M}}
      &\; \llb \phi \sto \psi \rrb^{\mf{M}}
      &= \llb \phi \rrb^{\mf{M}} \undsto \llb \psi \rrb^{\mf{M}}
  \end{align*}
  If $x \in \llb \phi \rrb^{\mf{M}}$ we say that $x$ \emph{satisfies} $\phi$,
  and write $\mf{M}, x \Vdash \phi$.
  A model $\mf{M}$ satisfies $\phi$ if $\llb \phi \rrb^{\mf{M}} = X$,
  notation: $\mf{M} \Vdash \phi$.
  A frame $\mf{X}$ is said to satisfy $\phi$
  if every model based on it satisfies $\phi$, in which case we write
  $\mf{X} \Vdash \phi$.
  If $\Gamma$ is a set of $\lan{L}_{\sto}$-formulae then we write
  $\mf{X} \Vdash \Gamma$ if $\mf{X} \Vdash \phi$ for all $\phi \in \Gamma$.
\end{definition}

  The class of $\sto$-models can be extended to a category with the following
  notion of morphism.

\begin{definition}
  A \emph{$\sto$-model morphism} from $(\mf{X}, V)$ to $(\mf{X}', V')$ is
  a $\sto$-frame morphism $f:\mf{X}~\to~\mf{X}'$ that additionally satisfies
  $V = f^{-1} \circ V'$.
  We write $\sMod$ for the category of $\sto$-models and their morphisms.
\end{definition}

  A routine induction on the structure of $\phi$ shows that:

\begin{proposition}\label{prop:mor-pres-truth}
  Let $f : \mf{M} \to \mf{M}'$ be a $\sto$-model morphism.
  Then for all $x \in X$ and $\phi \in \lan{L}_{\sto}$ we have
  $$
    \mf{M}, x \Vdash \phi \jff \mf{M}', f(x) \Vdash \phi.
  $$
\end{proposition}

\begin{corollary}\label{cor:mor-pres-truth}
  Let $f : \mf{X} \to \mf{X}'$ be a surjective $\sto$-frame morphism.
  Then $\mf{X} \Vdash \phi$ implies $\mf{X}' \Vdash \phi$.
\end{corollary}
\begin{proof}
  If $V'$ is any valuation for $\mf{X}'$, then $V := f^{-1} \circ V'$
  is a valuation for $\mf{X}$ such that
  $f : (\mf{X}, V) \to (\mf{X}', V')$ is a $\sto$-model morphism.
  It then follows from Prop.~\ref{prop:mor-pres-truth} that
  $(\mf{X}', V') \Vdash \phi$ because $(\mf{X}, V) \Vdash \phi$ and $f$ is surjective.
\end{proof}

  It is known (see e.g.~\cite[Thm.~6.4a]{LitVis18} and references therein)
  that $\sto$-frames form a sound and complete semantics for $\iA$.
  We now prove correspondence results for the axioms of $\lna{IELE}$
  and use this to give a sound and complete semantics for $\lna{IELE}$. As an example application, we substantiate the claim that $\lna{IELE}$ does not satisfy
  \ref{ax:Box}, \ref{ax:Hug}, and the converse of \ref{ax:IR}.

\begin{proposition}\label{prop:corr}
  Let $\mf{X} = (X, \preceq, \sqsubset)$ be a $\sto$-frame.
  \begin{enumerate}[label=\textup{\arabic*)}]
    \item \label{it:prop:corr-1}
          $\mf{X}$ satisfies \ref{ax:Sa} iff $x \sqsubset y$ implies
          $x \preceq y$ for all $x, y \in X$.
    \item \label{it:prop:corr-2}
          $\mf{X}$ satisfies \ref{ax:IR} iff for all $x \in X$ there exists
          $y \in X$ such that $x \preceq y$ and $x \sqsubset y$.
  \end{enumerate}
\end{proposition}
\begin{proof}
  Item \ref{it:prop:corr-1} follows from Lem.~4.10 and Thm.~6.4(c) in \cite{LitVis18}.

  For the second item, suppose there exists $x \in X$ such that no $\sqsubset$-successor of $x$
  is also a $\preceq$-successor. Then we can set $V(p) = {\uparrow}_{\preceq}x$
  and $V(q) = \emptyset$, and an easy verification shows that
  $x \Vdash \phi \sto \psi$ while $x \not\Vdash \phi \to \neg\neg\psi$.
  Conversely, let $\mf{X}$ be such that every $x \in X$ has a
  $\sqsubset$-successor that is also a $\preceq$-successor.
  Suppose $x \Vdash \phi \sto \psi$.
  We aim to prove that $x$ satisfies $\phi \to \neg\neg\psi$.
  So let $y$ be such that $x \preceq y$ and $y \Vdash \phi$.
  Then in order to prove that $y \Vdash \neg\neg\psi$, we need to
  show that $y \preceq z$ implies $z \not\Vdash \psi \to \bot$.
  By assumption there exists $w \in X$ such that
  $z \sqsubset w$ and $z \preceq w$.
  Since $x \preceq z$ and $x \Vdash \phi \sto \psi$ we must have
  $z \Vdash \phi \sto \psi$. Similarly, as $y \preceq w$ and $y \Vdash \phi$
  we find $w \Vdash \phi$, and because $z \sqsubset w$ this means $w \Vdash \psi$.
  Therefore $z \not\Vdash \psi \to \bot$.
  This proves that $y \Vdash \neg\neg\psi$ and hence
  $x \Vdash (\phi \to \neg\neg\psi)$, as desired.
\end{proof}

  In the presence  of the strength axiom, validity of \ref{ax:IR} is guaranteed
  by the requirement that every $x \in X$ has a $\sqsubset$-successor.
  Therefore a $\sto$-frame $(X, \preceq, \sqsubset)$ is a frame for $\lna{IELE}$
  if and only if $x \sqsubset y$ implies $x \preceq y$
  and every state has a $\sqsubset$-successor.
  We note that these are the same frames as used for the intuitionistic
  epistemic logic $\lna{IEL}$, cf.~\cite[Def.~4.2]{ArtPro16}.

\begin{proposition}\label{prop:IELE-form2}
  The axioms \ref{ax:Box}, \ref{ax:Hug}, and the converse of \ref{ax:IR}
  are not derivable from $\lna{IELE}$.
\end{proposition}
\tlnt{Make frame conditions explicit?}

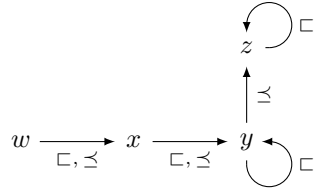
\begin{figure}
  \centering
    \begin{tikzpicture}
      \draw (0,0) node {$w$}
            (1.5,0) node {$x$}
            (3,0) node {$y$}
            (3,1.25) node {$z$};
      \draw[-latex] (.25,0) to node[anchor=north]{\scriptsize{${\sqsubset}, {\preceq}$}} (1.25,0);
      \draw[-latex] (1.75,0) to node[anchor=north]{\scriptsize{${\sqsubset}, {\preceq}$}} (2.75,0);
      \draw[-latex] (3,.25) to node[anchor=west]{\scriptsize{$\preceq$}} (3,1);
      \draw[] (3.0,-.3) arc(-180:90:.3);
      \draw[-latex] (3,-.25) -- (3,-.3) (3.3,0) -- (3.2,0);
      \draw (3.8,-.3) node {\scriptsize{$\sqsubset$}};
      \draw (3.3,1.25) arc(-90:180:.3);
      \draw[-latex] (3.25,1.25) -- (3.3,1.25) (3,1.55) -- (3,1.45);
      \draw (3.8,1.55) node {\scriptsize{$\sqsubset$}};
    \end{tikzpicture}
    \caption{The frame from Prop.~\ref{prop:IELE-form2}.}
    \label{fig:prop:IELE-form2}
\end{figure}

\jim{not in final LICS, yes in technical report}
\begin{proof}
  Consider the frame $\mf{X} = (X, \preceq, \sqsubset)$ where
  $$X = \{ w, x, y, z \}, \quad {\sqsubset} = \{ (w, x), (x, y), (y, y), (z, z) \},$$%
  and $\preceq$ is given by the reflexive and transitive closure of
  $w \preceq x \preceq y \preceq z$, see Fig.~\ref{fig:prop:IELE-form2}.
  This is a frame for $\lna{IELE}$ because the relation $\sqsubset$
  is contained in $\preceq$, and every state has a $\sqsubset$-successor.
  It can easily be verified that this does not satisfy the frame conditions
  corresponding to \ref{ax:Box} and \ref{ax:Hug} given in \cite[Fig.~6.2]{LitVis18},
  so that $\mf{X} \not\Vdash \axref{ax:Box}$ and $\mf{X} \not\Vdash \axref{ax:Hug}$.
  It follows that \ref{ax:Box} and \ref{ax:Hug} are not derivable in $\lna{IELE}$.
  
  Finally, we show that over $\mf{X}$,  
  $(\phi \to \neg\neg\psi)$ does not imply $(\phi \sto \psi)$.
  Consider the valuation given by $V(p) = X$ and $V(q) = \{ z \}$.
  Then by construction $x \not\Vdash p \sto q$.
  On the other hand, it follows from $z \Vdash q$ that $t \not\Vdash \neg q$ for
  all states $t \in X$, so that $t \Vdash \neg\neg q$ for all $t \in X$.
  Therefore $x \Vdash p \to \neg\neg q$, and hence $x$ itself witnesses
  $x \not\Vdash (p \to \neg\neg q) \to (p \sto q)$.
\end{proof}

\subsection{General and Descriptive Frames}\label{subsec:general}

  \noindent
  Next we define general $\sto$-frames. These will then be used to obtain a
  duality result for $\HLAs$ (Thm.~\ref{thm:duality}).
  Moreover, we will make extensive use of general
  frames in \rfse{embedding} below, where we embed $\lan{L}_{\sto}$
  into bimodal classical logic.

\begin{definition}
  A \emph{general $\sto$-frame} is a tuple $(X, \preceq, \sqsubset, P)$
  such that $(X, \preceq, \sqsubset)$ is a $\sto$-frame and
  $P \subseteq \fun{up}(X, \preceq)$ is a collection of upsets
  containing $X$ and $\emptyset$, and closed under
  $\cap, \cup, \undto$ and $\undsto$.
  It is called \emph{descriptive} if additionally it is
  \begin{itemize}
    \item \emph{Compact:} 
          For every $A \subseteq P$ and $B \subseteq \{ X \setminus a \mid a \in P \}$,
          if $A \cup B$ has the f.i.p.~then
          $\bigcap (A \cup B) \neq \emptyset$;
    \item \emph{$\preceq$-Refined:}
          For all $x, y \in X$, if $x \not\preceq y$ then there exists $a \in P$
          such that $x \in a$ and $y \notin a$;
    \item \emph{$\sqsubset$-Refined:}
          For all $x, y \in X$, if $x \not\sqsubset y$ then there exist $a, b \in P$
          such that $x \in a \undsto b$ and $y \in a$ and $y \notin b$.
  \end{itemize}
  If $\mf{G} = (X, \preceq, \sqsubset, P)$ is a geneneral $\sto$-frame,
  we write $\kappa\mf{G}$ for the underlying $\sto$-frame $(X, \preceq, \sqsubset)$.
\end{definition}

  Observe that the reduct $(X, \preceq, P)$ of a general $\sto$-frame
  $\mf{G} = (X, \preceq, \sqsubset, P)$ is a general intuitionistic Kripke frame.
  Therefore the set $P$ of admissibles forms a sub-Heyting algebra of $\fun{up}(X, \preceq)$.
%
  Moreover, if $\mf{G}$ is descriptive then $(X, \preceq, P)$ is a
  descriptive intuitionistic Kripke frame (see e.g.~\cite[\S\ 8.4]{ChaZak97}).
  Therefore we may alternatively define a descriptive $\sto$-frame as a tuple
  $(X, \preceq, \sqsubset, P)$ such that (i) $(X, \preceq, P)$ is a descriptive
  intuitionistic Kripke frame, (ii) $P$ is closed under $\undsto$,
  and (iii) $\sqsubset$-refinedness is satisfied.
  
  Since $P$ is closed under $\undsto$, we can view $(P, \undsto)$
  as a sub-algebra of $(\kappa\mf{G})^+$ (cf.~Lem.~\ref{lem:fun-plus-obj}).
  In particular this implies that $(P, \undsto)$ is an \HL-algebra,
  and we denote it by $\mf{G}^* = (P, \undsto)$.
  
  We now define morphisms between general $\sto$-frames.

\begin{definition}
  A \emph{general $\sto$-frame morphism} between $(X, \preceq, \sqsubset, P)$ and
  $(X', \preceq', \sqsubset', P')$ is a $\sto$-frame morphism
  $f : (X, \preceq, \sqsubset) \to (X', \preceq', \sqsubset')$ with 
  $f^{-1}(a') \in P$ for all $a' \in P'$.
  Let $\GFrm$ be the category of general $\sto$-frames and
  mor\-phisms and $\DFrm$ its full subcategory of descriptive $\sto$-frames.
\end{definition}

\begin{remark}\label{rem:dfrm-esa}
  It is well known that descriptive intuitionistic Kripke frames can be
  viewed as topological spaces with an extra relation, called
  Esakia spaces \cite{Esa74,Esa19}.
  As descriptive $\sto$-frames are based on descriptive intuitionistic
  Kripke frames, this adapts accordingly. 
  Define a \emph{strict implication space} to be a tuple $(X, \preceq, \sqsubset, \tau)$
  such that $(X, \preceq, \tau)$ is an Esakia space and $\sqsubset$ a
  binary relation on $X$ such that
  \begin{itemize}
    \item $x \preceq y \sqsubset z$ implies $x \sqsubset z$ for all $x, y, z \in X$;
    \item ${\downarrow}_{\sqsubset}a = \{ x \in X \mid x \sqsubset y \text{ for some }y \in a \}$
          is clopen for every clopen $a \subseteq X$;
    \item ${\uparrow}_{\sqsubset}x = \{ y \in X \mid x \sqsubset y \}$ is closed in $(X, \tau)$
          for all $x \in X$.
  \end{itemize}
  These constitute the category $\cat{SIS}$, whose morphisms are continuous
  morphisms that are bounded with respect to both relations.
  An easy verification
  shows that
  $$
    \DFrm \cong \cat{SIS}.
  $$
  For details, see Appendix \ref{app:descr-esa}.
\end{remark}

  If $f : \mf{G} \to \mf{G}'$ is a general $\sto$-frame morphisms,
  then we define $f^* := f^{-1} : (\mf{G}')^* \to \mf{G}^*$.
  This assignment is well defined by the definition of a general $\sto$-frame
  morphism, and it is an \HL-algebra morphism as a consequence of
  Prop.~\ref{prop:fun-plus}.
  We have:

\begin{proposition}
  The assignment $(\cdot)^* : \GFrm \to \HLAs$ defines a contravariant
  functor.
\end{proposition}

  A general $\sto$-frame $\mf{G} = (X, \preceq, \sqsubset, P)$
  can be turned into a \emph{general $\sto$-model} by
  endowing it with an \emph{admissible valuation}, that is, a map
  $V : \Prop \to P$. The interpretation of $\lan{L}_{\sto}$-formulae in
  $(\mf{G}, V)$ is defined as in the underlying $\sto$-model $(\kappa\mf{G}, V)$.
  We write $\mf{G} \Vdash \phi$ if $\phi$ is satisfied in every general
  $\sto$-model based on $\mf{G}$.
  If $\Gamma$ is a set of $\lan{L}_{\sto}$-formula then we define
  $\mf{G} \Vdash \Gamma$ as expected.
  Similar to Cor.~\ref{cor:mor-pres-truth} one can prove:
  
\begin{proposition}\label{prop:gen-mor-pres-truth}
  Let $f : \mf{G} \to \mf{G}'$ be a surjective general frame morphism.
  Then $\mf{G} \Vdash \phi$ implies $\mf{G}' \Vdash \phi$.
\end{proposition}

  Since valuations for $\mf{G}$ are in particular valuations for the underlying
  $\sto$-frame $\kappa\mf{G}$, validity of a formula $\phi$ in $\kappa\mf{G}$
  implies validity of $\phi$ in $\mf{G}$. The converse, however, need not
  be true. If the converse holds for the class of descriptive $\sto$-frames,
  then we call $\phi$ \emph{canonical}. Similarly, for logics:
  
\begin{definition}
  A logic $\Lambda$ is called \emph{canonical} if $\mf{G} \Vdash \Lambda$
  implies $\kappa\mf{G} \Vdash \Lambda$ for all descriptive $\sto$-frames $\mf{G}$.
\end{definition}

  This notion of canonicity is sometimes called \emph{d-persistence},
  as is the case in \cite{WolZak97,WolZak98}.

\begin{example}
  An example of a canonical $\iA$-logic is $\iA \oplus \axref{ax:Sa}$.
  To see this, note that by Prop.~\ref{prop:corr}\eqref{it:prop:corr-1}
  it suffices to prove that for every descriptive frame $\mf{G} = (X, \preceq, \sqsubset, P)$
  satisfying \ref{ax:Sa} we have: if $x \sqsubset y$ then
  $x \preceq y$ for all $x, y \in X$.
  This, in turn, follows from an application of $\sto$-refinedness.
\end{example}

  In \S\S~\ref{sec:embedding} and \ref{sec:bimodal} below we will see how
  canonicity for $\iA$-logics follows from canonicity of classical bimodal logics.

\subsection{Duality}\label{subsec:duality}
  
  \noindent
  Recall that every Heyting algebra $A$ gives rise to a descriptive
  intuitionistic Kripke frame $(\fun{pf}A, \subseteq, \tilde{A})$.
  Here $\fun{pf}A$ denotes the set of prime filters of $A$, and
  $\tilde{A} = \{ \tilde{a} \mid a \in A \}$, where
  $\tilde{a} = \{ \mf{p} \in \fun{pf}A \mid a \in \mf{p} \}$.
  (For details see e.g.~\cite[\S 8.2]{ChaZak97}
  or \cite[\S\S\ 2.2 and 2.3]{Bez06}.)
  We generalise this to construct a general $\sto$-frame from an \HL-algebra.
  This extends to a functor $(\cdot)_* : \HLAs \to \GFrm$,
  which is shown to give rise to a duality $\HLAs \equiv^{\op} \DFrm$
  in Thm.~\ref{thm:duality} below.
  
  This duality can also be obtained from \cite[Thm.~4.15]{CelJan05}
  by restricting to the right-hand side to WH-spaces \cite[Def.~4.4]{CelJan05}
  whose underlying Priestly space is an
  Esakia space. 

\begin{definition}\label{def:alg-dual}
  Let $\alg{A} = (A, \sto)$ be an \HL-algebra and let
  $(\fun{pf}A, \subseteq, \tilde{A})$ be the descriptive intuitionistic Kripke
  frame dual to $A$.
  Then we define the general $\sto$-frame $\alg{A}_*$ dual to $\alg{A}$ to be
  $\alg{A}_* = (\fun{pf}A, \subseteq, \sqsubset, \tilde{A})$, where
  $$
    \mf{p} \sqsubset \mf{q}
      \quad\text{iff}\quad
      \forall a, b \in A ( a \sto b \in \mf{p}
      \text{ and } a \in \mf{q}
      \text{ implies } b \in \mf{q}).
  $$
\end{definition}

  While it is fairly easy to verify that $\subseteq$ and $\sqsubset$ satisfy
  the required coherence condition \eqref{eq:p-sto}, so that
  $(\fun{pf}A, \subseteq, \sqsubset)$ is a $\sto$-frame,
  it is not clear whether $\tilde{A}$ gives it a general frame structure.
  In particular, it is not obvious why $\tilde{A}$ should be closed under
  $\undsto$.
  In order to prove that it is, we make use of the fact that
  $\tilde{\cdot} : A \to \tilde{A}$ defines an isomorphism of Heyting algebras
  \cite[\S~8.4]{ChaZak97}.
  In Lem.~\ref{lem:key} we prove that that it is a $\iA$-morphism,
  so that consequently 
  $\tilde{a} \undsto \tilde{b} = \widetilde{a \sto b} \in \tilde{A}$.
%

\begin{lemma}\label{lem:key}
  Let $\alg{A} = (A, \sto)$ be an \HL-algebra, and let
  $\alg{A}_* = (\fun{pf}A, \subseteq, \sqsubset, \tilde{A})$ be defined as
  in Definition \ref{def:alg-dual}.
  Then for all $a, b \in A$ we have
  $$
    \widetilde{a \sto b} = \tilde{a} \undsto \tilde{b}.
  $$
\end{lemma}
\begin{proof}
  A proof can be found in Appendix \ref{subsec:proof-lem-key}.
\end{proof}
%

  As stated, this shows that $\alg{A}_*$ is a general $\sto$-frame.
  In fact:

\begin{proposition}\label{prop:iA-to-gfrm}
  If $\alg{A} = (A, \sto)$ is a \HL-algebra,
  then $\alg{A}_*$ is a descriptive $\sto$-frame.
\end{proposition}

  As a consequence, the map $\tilde{(\cdot)} : \alg{A} \to (\alg{A}_*)^*$
  is an isomorphism in $\HLAs$.
  Indeed, it is an isomorphism between the underlying
  Heyting algebras, so by Lem.~\ref{lem:key} a bijective \HL-algebra morphism,
  hence an isomorphism because $\HLAs$ is a variety of algebras.

  Every $\iA$-morphism $ h : \alg{A} \to \alg{A}'$ is in particular
  a Heyting homomorphism, wherefore
  $\fun{pf}h = h^{-1} : (\fun{pf}A', \subseteq, \tilde{A}') \to (\fun{pf}A, \subseteq, \tilde{A})$
  is a general intuitionistic Kripke frame morphism.
  It follows from the definition of $\sqsubset$ that $h^{-1}$ is also a $\sto$-frame
  morphism between the underlying $\sto$-frames, so that setting
  $h_* = h^{-1}$ yields a functor $(\cdot)_* : \HLAs \to \GFrm$.
  
  
\begin{lemma}\label{lem:nat-iso-iA}
  The iso $\tilde{(\cdot)} : \alg{A} \to (\alg{A}_*)^*$ is natural
  in $\alg{A}$.
\end{lemma}

%
%

  Since the functor $(\cdot)_*$ lands in $\DFrm$, we may view it as a functor
  $\HLAs \to \DFrm$.
  In the converse direction we shall be sloppy and write $(\cdot)^*$ for the
  restriction of $(\cdot)^* : \GFrm \to \HLAs$ to
  $\DFrm$. Then we get:
  
\begin{theorem}\label{thm:duality}
  The functors $(\cdot)^* : \DFrm \to \HLAs$
  and $(\cdot)_* : \HLAs \to \DFrm$ define a
  dual equivalence
  $
    \HLAs \equiv^{\op} \DFrm.
  $
\end{theorem}

  As a consequence of the duality of Thm.~\ref{thm:duality} we have:
  
\begin{theorem}
  Every $\iA$-logic is characterised by a class of descriptive $\sto$-frames.
\end{theorem}

\section{Embedding into Bimodal Classical Logic}\label{sec:embedding}

  \noindent
  The G\"{o}del-McKinsey-Tarski translation embeds intuitionistic logic
  into the modal classical logic $\lna{S4}$ by prefixing every subformula
  of an intuitionistic formula with $\Box$ \cite{God33,McKTar48}.
  This was extended to the class of all intermediate logics by
  Dummett and Lemmon \cite{DumLem59}.
  The structure of the lattice of intermediate logics was later investigated by
  Maksimova and Rybakov \cite{MakRyb74},
  Blok \cite{Blo76}
  and Esakia \cite{Esa79,Esa19}.
  Most notably, this led to what is now known as the Blok-Esakia theorem,
  which classifies all modal companions of an intermediate logic and,
  as a corollary, establishes an isomorphism between the lattice of intermediate
  logics and the lattice of $\lna{S4}$-logics that satisfy the Grzegorczyk axiom.
  Fischer Servi and Shehtman further generalised it to embed intuitionistic
  logic with a normal unary modality into the bimodal classical logic
  $\lna{S4} \otimes \lna{K}$ \cite{Fis77,She79}.
  This was then exploited by Wolter and Zakharyaschev  \cite{WolZak98,WolZak97}
  to transfer results on completeness, decidability, the finite model property, and tabularity
  between modal intuitionistic logic and bimodal classical logic.
  
  In this section we generalise the G\"{o}del-McKinsey-Tarski translation
  to an embedding $t$  of $\lan{L}_{\sto}$ into a bimodal classical language. 
  We briefly recall some facts, fix notation
  for bimodal classical logic and define our syntactic translation in \rfsse{gen-sfk-frm}.
  Thereafter, in \rfsse{trans-frm}, we define translations from
  general $\sto$-frames to general $\sfk$-frames and vice versa, and examine
  their properties.
  In \rfsse{modal-comp} we define and investigate modal companions
  of $\iA$-logics,
  and prove that modal companions reflect decidability, Kripke completeness
  and the finite model property. That is, in order to prove that
  an $\iA$-logic has one of these properties it suffices to find a modal
  companion that does.
%
  

\subsection{General $\sfk$-frames}\label{subsec:gen-sfk-frm}

  \noindent
  We denote 
  the fusion of $\sf{S4}$ (with modality $\Boxi$) and $\sf{K}$
  (with modality $\Boxm$) by $\sfk$. The subscript $\sf{i}$ indicates that 
  $\Boxi$ arises from embedding \textbf{i}ntuitionistic logic into $\sf{S4}$.
  The box with subscript $\sf{m}$ is an additional \textbf{m}odality used for the translation
  of $\sto$. 

\begin{definition}
  An \emph{$\sfk$-frame} is a tuple $(X, \Ri, \Rm)$ consisting of
  a set $X$, and a pre-order $\Ri$ and a binary relation $\Rm $ on $X$.
  A \emph{p-morphism} from $(X, \Ri, \Rm)$ to $(X', \Ri', \Rm')$ is a function
  $f : X \to X'$ that is bounded with respect to both relations.
  An $\sfk$-frame $(X, \Ri, \Rm )$ is called a \emph{bimodal Heyting-Lewis frame}
  or \emph{$\sflb$-frame} if it satisfies
  \begin{equation}\label{eq:sflb} 
    \Ri \circ \Rm  \subseteq \Rm 
  \end{equation}
  We write $\cat{S4K}$ and $\cat{S4BHL}$ for the categories of $\sfk$-frames
  and $\sflb$-frames, respectively, and p-morphisms.
\end{definition}

  Observe that every $\sto$-frame can be conceived of as an $\sflb$-frame.
  However, the latter are still slightly more general than $\sto$-frames because
  $\Ri$ is only required to be a pre-order, rather than a partial order.
  An easy verification shows that satisfaction of \eqref{eq:sflb} is equivalent 
  to validity of 
  \begin{enumerate}[\qquad 1 \;]
    \axiom{BHL} \label{ax:HL}
                     $\Boxm\phi \to \Boxi\Boxm\phi$
  \end{enumerate}
  
  Furthermore, standard  Sahlqvist-style results \tlnt{ref} entail that the normal bimodal logic $\sflb$ obtained by extending $\sfk$ with \ref{ax:HL} is in fact (strongly) complete with respect to $\sflb$ frames. 
  
  We now recall the definition of general $\sfk$-frames, for a textbook reference
  see e.g.~\cite[\S~5]{BlaRijVen01} or \cite[\S~8.1]{ChaZak97}.

\begin{definition}\label{def:gen-sfk-frm}
  A \emph{general $\sfk$-frame} is a tuple $(X, \Ri, \Rm , P)$ that consists of an
  $\sfk$-frame $(X, \Ri, \Rm)$ and a Boolean subalgebra $P \subseteq \fun{P}X$
  of the powerset Boolean algebra of $X$, such that $P$ is closed under
  \begin{align*}
    \ibx : \fun{P}X \to \fun{P}X &: a \mapsto \{ x \in X \mid x\Ri y \text{ implies } y \in a \}, \\
    \mbx : \fun{P}X \to \fun{P}X &: a \mapsto \{ x \in X \mid x\Rm y \text{ implies } y \in a \}.
  \end{align*}
  A \emph{general $\sfk$-frame morphism} $f : (X, \Ri, \Rm , P) \to (X', \Ri', \Rm ', P')$
  is a p-morphism between the underlying $\sfk$-frame such that
  $f^{-1}(a') \in P$ whenever $a' \in P'$.
  We write $\cat{G\hyphen S4K}$ for the category of general $\sfk$-frames
  and morphisms.
  
  General $\sflb$-frames, descriptive $\sfk$-frames, and descriptive
  $\sflb$-frames are defined as usual (see e.g.~\cite[Def.~5.59 and 5.65]{BlaRijVen01}), 
  and their categories are denoted by $\cat{G\hyphen S4BHL}$,
  $\cat{D\hyphen S4K}$ and $\cat{D\hyphen S4BHL}$, respectively.
\end{definition}


  Finally, we write $\lna{Grz_i}$ \hypertarget{ax:Grz} and $\lna{Grz_m}$ for the 
  Grzegorczyk axiom
  \begin{enumerate}
    \axiom{Grz} $\Box(\Box(p \to \Box p) \to p) \to p$
  \end{enumerate} written with $\Boxi$ (resp. $\Boxm$) in place of
  $\Box$.

  Now we are ready to define (two versions of) the GMT translation. 

\begin{definition}[The GMT translation for $\lan{L}_{\sto}$]
  Let $\lanBM$ be the language of classical bimodal logic with modalities
  $\Boxi$ and $\Boxm$,
  and define the translation $t$ of $\lan{L}_{\sto}$ into
  $\lanBM$ recursively by
  \begin{align*}
    t(p) = \Boxi p &\qquad
    t(\top) = \top \qquad
    t(\bot) = \bot \\
    t(\phi \wedge \psi) &= \Boxi(t(\phi) \wedge t(\psi)) \\
    t(\phi \vee \psi) &= \Boxi(t(\phi) \vee t(\psi)) \\
    t(\phi \to \psi) &= \Boxi(t(\phi) \to t(\psi)) \\
    t(\phi \sto \psi) &= \Boxi\Boxm(t(\phi) \to t(\psi))
  \end{align*}
\end{definition}

\noindent 
  Note that the translation of $\Box\phi = \top \sto \phi$
  is given by
  $t(\Box\phi) = \Boxi\Boxm(t(\top) \to t(\phi))$, which is $\sfk$-equivalent to $\Boxi\Boxm t(\phi)$,
  so $t$ extends the translation used by Wolter and Zakharyaschev in \cite{WolZak98}. Furthermore, observe that over $\sflb$, the $\sto$-clause can be given as $\Boxm(t(\phi) \to t(\psi))$, a fact that we will use tacitly in what follows.

\subsection{Translations of Frames}\label{subsec:trans-frm}


  \noindent
  Next, we define functors $\hsigma : \GFrm \to \cat{G\hyphen S4K}$
  and $\hrho : \cat{G\hyphen S4K} \to \GFrm$.
  We prove that the composition $\hrho\hsigma$ is naturally isomorphic
  to the identity, and that both $\hsigma$ and $\hrho$ preserve descriptiveness.
  The transformations $\hrho$ and $\hsigma$ are based on the maps $\rho$
  and $\sigma$ that translate between general intuitionistic Kripke frames
  and general $\sf{S4}$-frames \cite[\S\S~3.9 and 8.3]{ChaZak97}.
  We add a hat to distinguish them from the maps $\rho$ and $\sigma$ that
  translate between $\iA$-logics and $\sfk$-logics defined in \rfsse{modal-comp}.
%
  A similar construction was carried out in \cite{WolZak97,WolZak98},
  and we point out the differences when we encounter them.

  We have already seen that every $\sto$-frame can be conceived of as an
  $\sfk$-frame. This extends to general frames:
%

\begin{definition}\label{def:sigma}
  Given a general $\sto$-frame $\mf{G} = (X, \preceq, \sqsubset, P)$,
  let $\hsigma\mf{G}$ be the general $\sfk$-frame
  $$
    \hsigma\mf{G} = (X, \preceq, \sqsubset, \hsigma P),
  $$
  where $\hsigma P$ is the Boolean closure of $P$ in $\fun{P}X$ (the powerset of $X$).
  For a general $\sto$-frame morphism $f$, let
  $\hsigma f = f$.
\end{definition}

\begin{lemma}\label{lem:sigmbas}
  If $\mf{G}$ is a general $\sto$-frame, then $\hsigma\mf{G}$ is a general $\sflb$-frame.
  Moreover, for any $a \in \hsigma P$ we have $\mbx a \in P$.
\end{lemma}
\begin{proof}
  We need to show that $\hsigma P$ is closed under $\ibx$ and $\mbx$.
  The former follows from \cite[Lem.~8.32 and 8.33]{ChaZak97}.
  For the latter (and the ``moreover'' part),
  first deconstruct $a \in \hsigma P$ as
  $$
    a = (\neg b_1 \cup c_1) \cap \dots \cap (\neg b_n \cup c_n)
  $$
  where $b_i$'s and $c_i$'s are elements of the original $P$. Then
  $$
    \mbx a = \mbx(\neg b_1 \cup c_1) \cap \dots \cap \mbx(\neg b_n \cup c_n).
  $$ 
  For each $i$, $\mbx(\neg b_i \cup c_i)$ is just $b_i \undsto c_i \in P$.
  The conclusion then follows from the fact 
  $P$ is closed under $\undsto$ and $\cap$.
\end{proof}
\noindent
  A straightforward verification shows that
  $\hsigma f : \hsigma\mf{G} \to \hsigma\mf{G}'$ is a general $\sflb$-frame morphism,
  whenever $f : \mf{G} \to \mf{G}'$ is a general $\sto$-frame morphism,
  and that $\hsigma$ is functorial, so that:

\begin{proposition}
  The assignment $\hsigma$ from Definition \ref{def:sigma} defines a
  functor $\hsigma : \GFrm \to \cat{G\hyphen S4BHL}$,
  and hence also from $\GFrm$ to $\cat{G\hyphen S4K}$.
\end{proposition}

%
  Now let us define a functor in the converse direction.
  For a general $\sfk$-frame $\mf{F} = (X, \Ri, \Rm , P)$, let
  $\Rm^* = {\Ri} \circ {\Rm}$.
  (This differs from \cite[\S 2]{WolZak98}, where \ref{ax:Box} is enforced by
  setting $\Rm^* = {\Ri} \circ {\Rm} \circ {\Ri}$.)
  Then $\mf{F}^* = (X, \Ri, \Rm^*, P)$ is a general $\sflb$-frame.
  The only thing separating $(X, \Ri, \Rm^*)$ from a $\sto$-frame is the fact
  that $(X, \Ri)$ is allowed to be a pre-order. 
  To resolve this, we quotient out cycles in $X$:
  
  Let $\mf{F} = (X, \Ri, \Rm , P)$ be a general $\sflb$-frame.
  The relation $\sim$ on $X$, given by $x \sim y$ if $x \Ri y$ and $y \Ri x$
  is an equivalence relation on $X$, whose equivalence classes are called
  \emph{($\Ri$-)clusters}.
  Let $\eqc{X}$ be the set of clusters of $(X, \Ri)$, and write $\eqc{x}$ for
  the cluster containing $x \in X$.
  Then $\Ri$ defines a partial order on $\eqc{X}$, which we denote by
  $\eqc{\Ri}$.
  Furthermore, define $\eqc{\Rm}$ by
  $$
    \eqc{x}\eqc{\Rm }\eqc{y} \jff x\Rm  y' \text{ for some } y' \sim y.
  $$
  Since  $x \sim x'$ implies $x\Ri x'$ and ${\Ri} \circ {\Rm} = {\Rm}$
  this does not depend on the choice of representative of $\eqc{x}$.
  
  For a set $a \subseteq X$ let $\eqc{a} = \{ \eqc{x} \mid x \in a \}$,
  and define
  $$
    \eqc{P} = \{ \eqc{a} \subseteq \eqc{X} \mid {\textstyle\bigcup} \eqc{a} \in P \}.
  $$
  Then an easy verification shows that
  $$\eqc{\mf{F}} = (\eqc{X}, \eqc{\Ri}, \eqc{\Rm}, \eqc{P})$$ is a general $\sflb$-frame.
  On passing, note that $\eqc{\mf{F}}$ is a general frame morphic image of
  $\mf{F}$ whenever $\mf{F}$ satisfies \ref{ax:HL}:

\begin{lemma}\label{lem:F-to-eqc-F}
  If $\mf{F}$ is a general $\sflb$-frame 
  then the assignment $x \mapsto \eqc{x}$ defines a general frame morphism
  $\mf{F} \to \eqc{\mf{F}}$.
\end{lemma}
\begin{proof}
  This follows immediately from the construction of $\eqc{\mf{F}}$.
\end{proof}

  Finally, we construct a general $\sto$-frame. Let
  $$
    \hrho P := \{ \ibx \eqc{a} \mid \eqc{a} \in \eqc{P} \}.
  $$
  It follows from \cite[\S 8.3]{ChaZak97} that $\hrho P$ is closed under
  $\cap, \cup$ and $\undto$. Moreover, $\hrho P$ is closed under $\undsto$ because
  $$
    \ibx a \undsto \ibx b
      = \mbx(\neg\ibx a \cup \ibx b)
      = \ibx\mbx(\neg\ibx a \cup \ibx b) \in \hrho P.
  $$
  Therefore we may define:
  
\begin{definition}
  For a general $\sfk$ frame $\mf{F} = (X, \Ri, \Rm , P)$ define
  $\hrho\mf{F} \in \GFrm$ by
  $$
    \hrho\mf{F} = (\eqc{X}, \eqc{\Ri}, \eqc{\Rm ^*}, \hrho{P}).
  $$
  For a morphism $f : \mf{F} \to \mf{F}'$ in $\cat{G\hyphen S4K}$ 
  define $\hrho f : \hrho\mf{F} \to \hrho\mf{F}'$ by 
  $\hrho f(\eqc{x}) = \eqc{f(x)}$.
\end{definition}

\begin{proposition}\label{prop:rho-fun}
  $\hrho : \cat{G\hyphen S4K} \to \GFrm$ is a functor.
\end{proposition}
\begin{proof}
  We have already seen that $\rho$ is well defined on objects.
  To see that the same goes for morphisms, 
  let $f : (X, \Ri, \Rm , P) \to (X', \Ri', \Rm ', P')$ be a morphism in 
  $\cat{G\hyphen S4K}$ and observe that
  monotonicity of $f$ proves that that $x \sim x'$ implies $f(x) \sim f(x')$,
  so that the definition of $\rho f$ does not depend on the choice of
  representative of $\eqc{x}$. Boundedness of $\rho f$ with respect to both
  relations is an immediate consequence of the fact that $f$ is a p-morphism.
  Furthermore, to see that $\rho f$ is a general frame morphism, we need to prove
  that $(\rho f)^{-1}(\ibx \eqc{a'}) \in \hrho P$ for all $a' \in P'$.
  This follows from the fact that
  $(\rho f)^{-1}(\ibx \eqc{a'}) = \ibx \eqc{f^{-1}(a')}$.

  Functoriality of $\rho$ is straightforward.
\end{proof}

\begin{remark}\label{rem:rho-sigma-full}
  We can also view $\hrho$ and $\hsigma$ as acting on non-general frames
  by viewing a $\sto$-frame as a general $\sto$-frame where every
  upset is admissible, and similar for $\sfk$-frames.
  This observation will be used in the proof of Thm.~\ref{thm:reflect}.
\end{remark}

  Not surprisingly, applying first $\hsigma$ and then $\hrho$ to a general
  $\sto$-frame yields an isomorphic frame.
  In fact, we can prove that the composition $\hrho\hsigma$ is naturally
  isomorphic to the identity functor on $\GFrm$.
  We will use this fact to prove facts about modal companions of $\iA$-logics
  in Thm.~\ref{thm:embeddable}.

\begin{proposition}\label{prop:rhosigma} 
  We have a natural iso $\hrho\hsigma \cong \fun{id}_{\GFrm}$.
\end{proposition}
\begin{proof}
  Let $\mf{G} = (X, \preceq, \sqsubset, P)$ be a general $\sto$-frame.
  By definition of $\rho$ and $\sigma$ we have
  $\hrho\hsigma\mf{G} = (X, \preceq, \sqsubset, \hrho\hsigma P)$,
  so for the isomorphism on objects we only have to show that $P = \hrho\hsigma P$.
  Since the definition of $\hrho\hsigma P$ is as in \cite[\S~8.3]{ChaZak97},
  this follows from Thm.~8.34 in {\it op.~\!cit.}
  Naturality of the isomorphism follows from the fact that $\hrho\hsigma f = f$.
\end{proof}

  As in \cite[Thm.~8.53]{ChaZak97}, we can prove that $\hsigma$ and $\hrho$
  preserve descriptiveness. This also extends \cite[Prop.~7]{WolZak98},
  but requires a more complicated proof.
  
  Descriptive frames are closely connected to Stone spaces, and we will make use
  of this perspective in the proof of Prop.~\ref{prop:pres-descr}.
  In particular, if $\mf{F} = (X, \Ri, \Rm, P)$ is a descriptive $\sfk$-frame
  and $\tau_A$ is topology on $X$ generated by (clopen) base $P$, then
  $(X, \tau_A)$ is a Stone space, hence compact, and 
  $\Ri[x] = \{ y \in X \mid x \Ri y \}$ and $\Rm[x] = \{ y \in X \mid x \Rm y \}$
  are closed in $\tau_A$ for all $x \in X$. See also \cite{KupKurVen04}.
  
\begin{proposition}\label{prop:pres-descr}
  $\hrho$ and $\hsigma$ preserve descriptiveness.
\end{proposition}
\begin{proof}[Proof of Proposition \ref{prop:pres-descr}]
  Suppose $\mf{G} = (X, \preceq, \sqsubset, P)$ is a descriptive $\sto$-frame.
  Then the fact that it is differentiated and compact proves that
  $\hsigma\mf{G}$ is differentiated and compact.
  Tightness follows from the fact that $\mf{G}$ is a descriptive
  $\sto$-frame.
  
  Conversely, suppose $\mf{F} = (X, \Ri, \Rm, P)$ is a descriptive $\sfk$-frame.
  Then $(\eqc{X}, \eqc{\Ri}, \hrho P)$ is a descriptive intuitionistic Kripke
  frame by \cite[Thm.~8.53]{ChaZak97}.
%
  So we only have to show that
  \begin{equation}\label{eq:descr}
  \begin{split}
    \eqc{x}\eqc{\Rm ^*}\eqc{y}
      \jff &\text{for all } \ibx{\eqc{a}}, \ibx \eqc{b} \in \hrho P, \\
      &\text{ if } \eqc{x} \in \ibx\eqc{a} \undsto \ibx\eqc{b} \\
      &\text{ and } \eqc{y} \in \ibx\eqc{a} \\
      &\text{ then } \eqc{y} \in \ibx\eqc{b}
  \end{split}
  \end{equation}
  The direction from left to right follows immediately from the definition of $\sto$,
  so we focus on the converse. 
  
  Suppose $\neg(\eqc{x}\eqc{\Rm^*}\eqc{y})$.
  Then $\neg(x\Rm^*y')$ for all $y' \sim y$.
  We will construct $a, b \in P$ that are up-closed under $\Ri$
  such that $y \in a$ and $y \notin b$ and $\Rm^*[x] \cap a \subseteq b$.
  (Here $\Rm^*[x] = \{ z \in X \mid x\Rm^*z \}$.)
  Then, since $a$ and $b$ are up-closed under $\Ri$, we have $\bigcup \eqc{a} = a$
  and $\bigcup \eqc{b} = b$, so that $\ibx\eqc{a}, \ibx\eqc{b} \in \hrho P$.
  Moreover, by construction $\ibx\eqc{a}$ and $\ibx\eqc{b}$ are such that
  $\eqc{x} \in \ibx\eqc{a} \undsto \ibx\eqc{b}$ and
  $\eqc{y} \in \ibx\eqc{a}$, while $\eqc{y} \notin \ibx\eqc{b}$.
  Therefore, they witness the right-to-left direction in \eqref{eq:descr}.
  
  We view $\mf{F} = (X, \Ri, \Rm, P)$ as a Stone space $(X, \tau_P)$ with
  point-closed relations $\Ri$ and $\Rm$ (cf.~\rfsse{gen-sfk-frm}).
  This allows us to use topological in the construction of $a$ and $b$.

  So suppose $\neg(x \Rm^* y')$ for all $y' \sim y$.
  Then ${\uparrow}_{\sf{m}}x = \{ z \in X \mid x\Rm z \}$
  and ${\downarrow}_{\sf{i}}y = \{ z \in X \mid z \Ri y \}$
  are closed, and hence their intersection
  $C = {\uparrow}_{\sf{m}}x \cap
  {\downarrow}_{\sf{i}}y$ is closed in $(X, \tau_P)$.
  Furthermore, by assumption $C$ does not contain any $y' \sim y$.
  Therefore none of the elements in $C$ lies above $y$ in the $\Ri$-ordering.
  We claim that we can find a clopen $\Ri$-upset $a$ containing $y$ and disjoint
  from ${\uparrow}_{\sf{m}}x \cap {\downarrow}_{\sf{i}}y$.
  To see this, note that since ${\uparrow}_{\sf{i}}y$ is closed we have
  ${\uparrow}_{\sf{i}}y = \bigcap \{ a \in P \mid {\uparrow}_{\sf{i}}y \subseteq a \}$.
  Therefore
  \begin{align*}
    {\uparrow}_{\sf{i}}y
      = \ibx ({\uparrow}_{\sf{i}}y)
      &= \ibx \Big(\bigcap \{ a \in P \mid {\uparrow}_{\sf{i}}y \subseteq a \}\Big) \\
      &= \bigcap \big\{ \ibx a \in P \mid {\uparrow}_{\sf{i}}y \subseteq a \big\}.
  \end{align*}
  A compactness argument using compactness of $C$ and the fact that
  $C$ and ${\uparrow}_{\sf{i}}y$ are disjoint now yields
  a clopen set $a := \ibx a_1 \cap \cdots \cap \ibx a_n \in P$
  (hence up-closed under $\Ri$)
  containing ${\uparrow}_{\sf{i}}y$ and disjoint from $C$.
  
  Similarly, using the fact that ${\uparrow}_{\sf{m}}x \cap a$ is closed
  and does not contain $y$,
  we can find a clopen $\Ri$-upset $b$ containing ${\uparrow}_{\sf{m}}x \cap a$
  such that $y \notin b$.
  Then we have $y \in \ibx a = a$, $y \notin \ibx b = b$ and
  ${\uparrow}_{\sf{m}}x \cap a \subseteq b$,
  as desired. 
\end{proof}

\subsection{Modal Companions}\label{subsec:modal-comp}

  \noindent
  Recall that $t : \lan{L}_{\sto} \to \lanBM$ denotes the extension of the
  G\"{o}del-McKinsey-Tarski translation of $\lan{L}_{\sto}$ into
  classical bimodal logic.
  The translation $t$ and functor $\hrho$  
  are related as follows.

\begin{lemma}\label{lem:translation-pres}
  Let $\mf{F} \in \cat{G\hyphen S4K}$ and $\phi \in \lan{L}_{\sto}$.
  Then
  $$
    \mf{F} \Vdash t(\phi) \jff \hrho\mf{F} \Vdash \phi.
  $$
\end{lemma}
\begin{proof}
  See Appendix \ref{subsec:translation-pres-proof}.
\end{proof}

  This lemma extends Lem.~8.28 in \cite{ChaZak97},
  and is the ``$\sto$-analogue'' of \cite[Lem.~5]{WolZak98}.
  It gives rise to the following (standard) notion of
  modal companions.

\begin{definition}
  Let $\Lambda$ be an extension of $\lna{iA}$ and
  $\Theta$ an extension of $\lna{S4K}$.
  If for all $\phi \in \mc{L}_{\sto}$ we have
  $$
    \phi \in \Lambda \jff t(\phi) \in \Theta
  $$
  then we say that $\Lambda$ is \emph{embedded} in $\Theta$,
  and $\Theta$ is an \emph{$\sfk$-companion} or \emph{modal companion} of $\Lambda$.
\end{definition}

  Analogously to \cite[Thm.~9.56]{ChaZak97},
  it follows from Lem.~\ref{lem:translation-pres} that
  for every \sfk-logic $\Theta$, the $\lna{iA}$-logic
  $$
    \rho\Theta = \{ \phi \in \mc{L}_{\sto} \mid t(\phi) \in \Theta \}
  $$
  is such that $\Theta$ is an $\sfk$-companion of $\rho\Theta$.
  Furthermore, it follows from Lem.~\ref{lem:translation-pres} that:

\begin{proposition}\label{prop:companion-char}
  If an $\sfk$-logic $\Theta$ is characterised by a class $\mc{C}$ of
  general $\sfk$-frames, then $\rho\Theta$ is characterised by the class
  $\hrho\mc{C} = \{ \hrho\mf{F} \mid \mf{F} \in \mc{C} \}$ of general
  $\sto$-frames.
\end{proposition}

  The proof of the following theorem resembles to proof of 
  \cite[Thm.~9]{WolZak98}, using Lem.~\ref{lem:translation-pres}
  and Prop.~\ref{prop:rhosigma}.

\begin{theorem}\label{thm:embeddable}
  Each $\iA$-logic $\Lambda = \iA \oplus \Gamma$ is embeddable by $t$ in any logic
  $\Theta$ in the interval
  $$
    [(\lna{S4} \otimes \lna{K}) \oplus t(\Gamma),
      (\Grzi \otimes \lna{K}) \oplus t(\Gamma) \oplus \text{\ref{ax:HL}}].
  $$
\end{theorem}

  We have seen how $\sfk$-logics give rise to $\iA$-logics.
  In the converse direction, guided by Thm.~\ref{thm:embeddable},
  we define:

\begin{definition}
  \jim{Not sure if we need this definition. Maybe delete it?}
  For an $\iA$-logic $\Lambda = \iA \oplus \Gamma$,
  let
  \begin{align*}
    \tau \Lambda &= (\lna{S4} \otimes \lna{K}) \oplus t(\Gamma) \oplus \axref{ax:HL} \\
    \sigma \Lambda &= (\Grzi \otimes \lna{K}) \oplus t(\Gamma) \oplus \axref{ax:HL}
  \end{align*}
\end{definition}

\begin{theorem}\label{thm:reflect}
  The map $\rho$ preserves decidability, Kripke completeness and the finite model
  property.
\end{theorem}
\begin{proof}
  Let $\Theta$ be an $\sfk$-logic.
  Since $\Theta$ is a modal companion of $\rho\Theta$ we have
  $\phi \in \rho\Theta$ if and only if $t(\phi) \in \Theta$.
  Therefore, deciding if $\psi \in \rho\Theta$ is equivalent to
  deciding if $t(\psi) \in \Theta$, and hence decidability of $\Theta$
  implies decidability of $\rho\Theta$.
  
  Next, we assume that $\Theta$ is Kripke complete.
  In order to prove that $\rho\Theta$ is Kripke complete as well, we show
  that for each $\psi \in \lan{L}_{\sto}$ such that $\psi \notin \rho\Theta$
  we can find a $\sto$-frame for $\rho\Theta$ that does not
  validate $\psi$.
  If $\psi \notin \rho\Theta$ then $t(\psi) \notin \Theta$.
  Since $\Theta$ is Kripke complete, we can find a \sfk-frame $\mf{F}$
  that validates $\Theta$ such that $\mf{F} \not\Vdash t(\psi)$.
  Recall from Remark~\ref{rem:rho-sigma-full} that we can view $\hrho$
  as a map taking \sfk-frames to $\sto$-frames.
  So we view $\hrho\mf{F}$ as a $\sto$-frame, and by
  Lemma~\ref{lem:translation-pres} it validates $\rho\Theta$ but not $\psi$.
  
  Lastly, we can prove preservation of the finite model property in the same
  way as preservation of Kripke completeness, using the observation that
  $\hrho$ sends finite frames to finite frames.
\end{proof}

\begin{examples} \label{ex:corr} \
  \begin{enumerate}
    \item 
          A modal companion 
          of  \ref{ax:Box} is given by the extension of $\sfk$ with
          $$
            \Boxi\Boxm\Boxi p \to \Boxm p.
          $$%
          This axiom is canonical and elementary via the SQEMA
          algorithm~\cite{Conradie:2006:ACC,Geo17} (or the usual Sahlqvist argument),
          which also yields strong completeness.
    \item \label{grztrick} The translation of  \ref{ax:P} from \rfse{HA} is
          $\sflb$-equivalent to $$\Boxm(\Boxi p \to \Boxi q) \to \Boxm\Boxm(\Boxi p \to \Boxi q).$$%
             Here we need to be somewhat creative. 
             Namely, we observe that by the results of this section, if $\tau\Lambda$ is canonical, then $\Lambda$ is strongly complete with respect to $\tau\Lambda$-frames, which are in addition partial orders. Next, one observes that over partial orders, the following rule is admissible: from $\phi(\Boxi p \to \Boxi q)$, derive $\phi(r)$, where $p$ and $q$ are fresh for $\phi(r)$
              The above translation can be verified to be canonical via the SQEMA algorithm \cite{Conradie:2006:ACC,Geo17}. In this way, we finally arrive at $\Rm$-transitivity as its (canonical) counterpart (see also \cite[\S\ 10]{LitVis19}).
     \item The
          translation of \ref{ax:IR} reads
          $$
            \Boxi(\Boxi\Boxm(\Boxi p \to \Boxi q) \to \Boxi(\Boxi p \to \Boxi\Diamondi\Boxi q)),
          $$%
          where $\Diamondi$ is short for $\neg\Boxi\neg$.
          As a consequence of Thm.~\ref{thm:embeddable} the logic
          $\tau\lna{IELE} = \sfk \oplus t(\axref{ax:Sb}) \oplus t(\axref{ax:IR})
          \oplus \axref{ax:HL}$
          is a modal companion of $\lna{IELE}$.
  \end{enumerate}
\end{examples}

\section{A Criterion for FMP and Decidability}\label{sec:bimodal}

\tadeusz{More and better comments where problems lie}

\noindent
Theorem~\ref{thm:reflect} allows transfer of decidability and the finite model property from a bimodal logic to its $\tto$-counterpart. Recall that for a finitely axiomatizable logic, the finite model property is a sufficient  criterion of decidability. Still, in order to use this theorem for a $\tto$-logic, we need a bimodal counterpart which enjoys these properties. Wolter and Zakharyaschev \cite{WolZak98} proposed a broad criterion based on  techniques  for unimodal (cofinal) subframe logics.
We begin by generalizing their criterion to $\sflb$-logics. They essentially relies on transitivity of \emph{both} relations: not only $\Ri$, but also $\Rm$. 
  Logically, the latter corresponds to validity of
  \begin{enumerate}
    \axiom{4_m} \label{ax:4m}
    $\Boxm p \to \Boxm\Boxm p$
  \end{enumerate}
  That is, we will be concerned with $(\lna{S4} \otimes \lna{K4})$-logics.
  In the presence  of the \ref{ax:Mix}-axiom, i.e.,
  \begin{enumerate}
    \axiom{Mix} \label{ax:Mix}
    $\Boxm p \to \Boxi\Boxm\Boxi p$
  \end{enumerate}
  the original criterion reads as
\begin{theorem}[\cite{WolZak98}, Thm.~17]
  Let $\Theta$ be a canonical subframe extension of $\sfk$ containing \ref{ax:Mix}. 
  If $\lna{S4} \oplus \Gamma \subseteq \lan{L}_{\sf{i}}$ is a $\Ri$-subframe logic,
  then $\Theta \oplus \Gamma$ has the finite model property.
\end{theorem}

  The goal of this section is to prove similar theorems that do not require
  \ref{ax:Mix}.
  First, recall that the subframe of $\mf{F} = (X, \Ri, \Rm, P)$ generated
  by $Y \subseteq X$ is the frame
  $\mf{F}_{\upharpoonright Y} = (Y, {\Ri}_{\upharpoonright Y}, {\Rm}_{\upharpoonright Y}, P_{\upharpoonright Y})$,
  where ${\Ri}_{\upharpoonright Y}$ and ${\Rm}_{\upharpoonright Y}$ are the
  restrictions of $\Ri$ and $\Rm$ to $Y$, and $P_{\upharpoonright Y} = \{ a \cap Y \mid a \in P \}$.
  It is called \emph{$\Rm$-cofinal} if for all $y \in Y$,
  $y \Rm z$ implies that there exists $y' \in Y$ such that $z \Rm y'$.
  An $\sfk$-logic $\Theta$ is called a \emph{($\Rm$-cofinal) subframe logic}
  if the collection of frames validating $\Theta$ is closed under forming
  ($\Rm$-cofinal) subframes.

  Besides, we make heavy use of the notions of \emph{$\Ri$- and $\Rm$-maximality} with respect to an equivalence
  relation generated by a formula $\phi$: 

\begin{definition}\label{def:maxim}
  Let $\mf{M} = (\mf{F}, V)$ be a model based on a general $\sfk$-frame
  $\mf{F} = (X, \Ri, \Rm, P)$.
  If $x \in X$ then we write ${\uparrow}_{\sf{i}}x = \{ y \in X \mid x \Ri y \}$
  for the upwards closure of $x$ and
  $$C_{\sf{i}}(x) = \{ y \in X \mid x \leq y \leq x \}$$ for the $\Ri$-cluster of $x$.
  The state $x$ is called \emph{$\Ri$-final}, and $C(x)$ is called an $\Ri$-final
  $\Ri$-cluster, if ${\uparrow}_{\sf{i}}x = C_{\sf{i}}(x)$.
  We similarly define ($\Rm$-final) $\Rm$-clusters.
  
  Let $\phi \in \lanBM$.
  We call $x, y \in X$ \emph{$\phi$-equivalent} in $\mf{M}$, and write
  $x \sim_{\phi} y$, if $x \Vdash \psi$ iff $y \Vdash \psi$ for all
  $\psi \in \Subf(\phi)$.
%
  A state $x$ is called $\Ri$-maximal in $\mf{M}$ (relative to $\sim_{\phi}$)
  if for any $x \neq y \in X$ such that $x \Ri y$ we have $x \not\sim_{\phi} y$.
  We similarly define $\Rm$-maximality.
\end{definition}

  If $\mf{M}$ is based on a descriptive
  frame, then for every state $x \in \mf{M}$ there exists an $\Ri$-maximal
  state $y$ such that $x \sim_{\phi} y$ and either $x = y$ or $x \Ri y$ \cite[Lem.~14]{WolZak98}.
  If $\Rm$ is transitive then the same holds for $\Rm$.
  

  The purpose of $\phi$-equivalence is showcased in the following lemma,
  that will prove useful in the subsequent theorem.
  
\begin{lemma}\label{lem:subsub-truth}
  Let $\mf{M} = (\mf{F}, V)$ be a model based on an $\sfk$-frame
  $\mf{F} = (X, \Ri, \Rm, P)$, and let $\phi \in \lanBM$.
  Suppose $Y \subseteq X$ is such that for $R \in \{ \Ri, \Rm \}$:
  if $y \in Y$ and $y R x$ (where $x \in X$),
          then there exists $y' \in Y$ such that $y R y'$ and $x \sim_{\phi} y'$.
  Then for all $y \in Y$ and $\psi \in \Subf(\phi)$ we have
  $$
    \mf{M}_{\upharpoonright Y}, y \Vdash \psi \jff \mf{M}, y \Vdash \psi.
  $$
\end{lemma}
\begin{proof}
  By induction on the structure of $\psi$.
  If $\psi = p, \top$ or $\bot$ then the statement is obvious,
  as is the case for $\phi = \psi_1 \wedge \psi_2$ and $\psi = \psi_1 \vee \psi_2$.
  
  If $\psi = \Boxi \psi_1$ and $\mf{M}, y \Vdash \Boxi\psi_1$, then
  clearly $\mf{M}_{\upharpoonright Y}, y \Vdash \Boxi\psi_1$, since every
  $({\Ri}_{\upharpoonright Y})$-successor of $y$ is also an $\Ri$-successor of $y$ in $\mf{M}$.
  Conversely, if $\mf{M}, y \not\Vdash \Boxi\psi_1$, then there is an 
  $\Ri$-successor $z$ such that $\mf{M}, z \not\Vdash \psi_1$.
  By assumption there exists $z' \in Y$ such that $y \Ri z'$ and $z \sim_{\phi} z'$.
  Consequently $\mf{M}, z' \not\Vdash \psi_1$, so by the induction
  hypothesis $\mf{M}_{\upharpoonright Y}, z' \not\Vdash \psi_1$
  and hence $\mf{M}_{\upharpoonright Y}, y \not\Vdash \Boxi \psi_1$.
  
  The case $\psi = \Boxm\psi_1$ is analogous.
\end{proof}

\newcommand{\cla}[1]{\textsf{#1}}
\newcommand{\mix}{\lna{mix}}

\begin{theorem}\label{thm:fmp-total}
  Suppose $\Theta$ is a canonical extension of $\lna{S4} \otimes \lna{K4}$
  containing \ref{ax:HL} that is closed under forming ($\Rm$-cofinal) subframes. Then:
  \begin{enumerate}[label=\textup{\arabic*)}]
    \item \label{it:thm:bm-1}
          $\Theta$ has the finite model property.
    \item \label{it:thm:bm-2}
          If moreover $\Theta$ contains the classical strength axiom
\begin{enumerate} 
 \axiom{S_c} \label{ax:Sc}
    $\Boxi p \to \Boxm p$.
 \end{enumerate}   
          then for any ($\Rm$-cofinal) subframe logic $\Gamma \subseteq \lan{L}_{\sf{m}}$,
          the logic $\Theta \oplus \Gamma$ has the finite model property.
  \end{enumerate}
\end{theorem}
\begin{proof}
  See Appendix \ref{app:proof-thm-fmp-total}.
\end{proof}

\begin{corollary}\
Let $\Lambda$ be a $\tto$-logic extending \loga{P}. 
\begin{enumerate}[label=\textup{\arabic*)}]
\item If its \sflb-counterparts include a canonical logic preserved by forming (cofinal) subframes, $\Lambda$ has the finite model property.
\item Furthermore,  if $\Lambda$ extends \loga{Sa} and its \sflb-counterparts include a logic  obtained by extending a canonical (cofinal) subframe logic  with a collection of  $\lan{L}_{\sf{m}}$-axioms preserved by $\Rm$-subframes, $\Lambda$ has the finite model property. 
\end{enumerate}
In either case, $\Lambda$ is decidable whenever finitely axiomatizable.
\end{corollary}

\begin{examples}\
\begin{enumerate}
\item The above theorem covers \loga{P} and \loga{Sa} themselves. As we have seen, their natural $\sflb$ couterparts are complete with respect to frames definable by universal first-order conditions, hence they are not only canonical, but also preserved under subframes.
\item It appears more challenging to use the second clause of the above corollary, as the GMT translation always includes some $\Boxi$ modalities. However, transformations similar to those used in Exm.~\ref{ex:corr} can handle, e.g., a syntactic variant of the strong L\"ob axiom
$$
(((p \sto q) \wedge p) \sto q) \to (p \sto q).
$$%
After the GMT translation, one application of the trick from Exm.~\ref{ex:corr}.\ref{grztrick} (plus some trivial book-keeping) yields L\"ob for $\Boxm$. This is a $\Rm$-subframe axiom.
\item Simlarly, when one considers $\lna{PLAA}$, its \ref{ax:App} axiom over \loga{Sa} translates to
$$
\Boxm((\Boxi p \wedge \Boxm(\Boxi p \to \Boxi q)) \to \Boxi q).
$$%
The usual currying trick yields
$$
\Boxm(\Boxm(\Boxi p \to \Boxi q)) \to (\Boxi p \to \Boxi q)),
$$%
which in one application of the  trick from  Exm.~\ref{ex:corr}.\ref{grztrick} produces
$
\Boxm(\Boxm r \to r)
$
and this is a $\Rm$-subframe axiom.
\end{enumerate}
\end{examples}

\noindent

\section{Conclusions and Future Work} \label{sec:conclusions}

\noindent
We have investigated the Heyting-Lewis (family of) logic(s) of strict implication.
We have described suitable categorical duality,
have provided a (truth preserving) translation into classical
bimodal logic, yielding results on canonicity, 
the finite model property, and decidability. 

Our study 
leaves many questions open.
Most obviously, we leave open:
\begin{question}
  Does the G\"odel-McKinsey-Tarski translation from Heyting-Lewis logic to
  bimodal classical logic give rise to an analogue of the Blok-Esakia theorem?
\end{question}
Other directions for further research include pushing the limits of the finite model
property results discussed in Section \rfse{bimodal},
and investigating
the Intuitionistic Logic of Entailment \lna{IELE}
proposed in \rfse{iele} in more detail. Apart from its
mathematical and philosophical aspects, or the question of full
axiomatization capturing all justifiable principles, we note
 the
Brouwer-Heyting-Kolmogorov motivation underlying the original
\lna{IEL}  \cite{ArtPro16}. Here,  one could relate it to type-theoretic calculi developed in the functional programming setting (cf. \rfse{haskell}). 
Finally, it would be of
interest to compare \lna{IELE} with other approaches to
intuitionistic epistemic logics \cite{Wil92,Pro12,KurPal13}.

On the more conceptual side, we are interested in viewing Heyting-Lewis
logic from a dialgebra perspective \cite{GroPat20}, in particular as
the semantics used in this paper naturally lends itself to a
dialgebraic treatment. 
  Dialgebraic methods can, for example, provide a generic
  expressivity and expressivity-somewhere-else results
  \cite[\S~7]{GroPat20} and Goldblatt-Thomason theorems
  \cite{Gro20b}.  It would be intriguing to see
  whether dialgebraic methods can provide an elegant connection
  (say, via a dialgebraic generalisation of the final coalgebra
  sequence) with, on the one hand, finitary completeness proofs
  developed in the preservativity community
  \cite{Iem01,IemhoffJZ05:igpl} and on the other hand, step-algebras
  and step-frames \cite{BezGhi14}.  



But what we see as the main future challenge is the development of semantics and decidability results of $\iAm$ logics not including \lna{Di}. We have seen that such logics arise under the functional programming interpretation (arrows without choice), and also under the arithmetical interpretation (e.g., preservativity without provable closure under q-realizability).\footnote{Especially in the functional programming world, it is common to consider fragments not involving disjunction, so the reader may ask if for such formulae at least, one cannot use Kripke completeness for an extension of a disjunction-free set of axioms with \lna{Di}; in other words, if extensions with \lna{Di} are conservative over the disjunction-free fragment. While the potential scope of such results is of independent interest, the general answer  is on the negative \cite[Lem. 4.6]{LitVis18}, \cite[\S\ 10]{LitVis19}.}
 Let us note here that $\iAm$ extends the minimal system \lna{ICK} of \emph{intuitionistic conditional logic} proposed recently by Weiss \cite{Wei19a,Wei19b,CiaLiu19}, which allows the use of an intuitionistic variant of Chellas frames in conditional logic \cite{Che80}.  However, we believe that the right approach is to generalize so-called Veltman semantics of classical \emph{interpretability logics} \cite{dJV90}, as classically preservativity  is the contrapose of arithmetical interpretability. Conceivable variants of generalized Veltman semantics would be equivalent to subclasses of Chellas-Weiss frames (much like Kripke frames can be seen as a limiting case of neighbourhood frames). Nevertheless, in the classical setting, (generalized) Veltman semantics has proved particularly suitable for decidability and complexity results \cite{dJV90,MikecPV17:igpl,MikecPV19:igpl}, allowing adaptations of standard modal techniques such as filtration \cite{PerkovV16}, so it does seem promising to work with more restrictive structures.  




\clearpage
\printbibliography

\clearpage
\appendix
\section{Appendix}

\noindent
  Here we collect details omitted from the main body of the text.
  Appendix \ref{app:descr-esa} elaborates on Remark \ref{rem:dfrm-esa}.
  Appendices
    \ref{subsec:proof-lem-key},
    \ref{subsec:translation-pres-proof}
    and 
    \ref{app:proof-thm-fmp-total}
  give proofs of
    Lem.~\ref{lem:key},
    \ref{lem:translation-pres}
    and Theorem \ref{thm:fmp-total}.
  These are given their own appendix because they require additional
  definitions and lemmas.
%

\subsection{Descriptive $\sto$-frames as spaces}\label{app:descr-esa}

  \noindent
  We work out the details of Remark \ref{rem:dfrm-esa}.
  Recall that a \emph{Stone space} is a zero-dimensional compact Hausdorff space.
  An \emph{Esakia space} \cite{Esa74} is a tuple $(X, \preceq, \tau)$ consisting of a
  Stone space $(X, \tau)$ and a relation $\preceq$ on $X$ that satisfies
  \begin{itemize}
    \item ${\uparrow}_{\preceq}x = \{ y \in X \mid x \preceq y \}$ is closed
          in $(X, \tau)$ for each $x \in X$;
    \item ${\downarrow}_{\preceq}a = \{ y \in X \mid y \preceq x \text{ for some } x \in a \}$
          is clopen in $(X, \tau)$ for all clopen $a \in \tau$.
  \end{itemize}
  Together with bounded continuous morphisms, Esakia spaces for the category
  $\cat{Esa}$. It is well known that $\cat{Esa}$ is isomorphic to the category
  of descriptive intuitionistic Kripke frames and general intuitionistic
  Kripke frame morphisms \cite{Esa74,Esa19}.
  
  Piggy-backing on this, we develop a topological perspective of descriptive
  $\sto$-frames.

\begin{definition}
  A \emph{strict implication space} is a tuple $(X, \preceq, \sqsubset, \tau)$
  such that $(X, \preceq, \tau)$ is an Esakia space and $\sqsubset$ is a
  binary relation on $X$ such that
  \begin{enumerate}[\qquad 1]
  \renewcommand{\theenumi}{($\mathsf{S_{\arabic{enumi}}}$)}
    \item \label{it:sis-1}
          $x \preceq y \sqsubset z$ implies $x \sqsubset z$ for all $x, y, z \in X$;
    \item \label{it:sis-3}
          ${\uparrow}_{\sqsubset}x = \{ y \in X \mid x \sqsubset y \}$ is closed in $(X, \tau)$
          for all $x \in X$;
    \item \label{it:sis-2}
          ${\downarrow}_{\sqsubset}a = \{ x \in X \mid x \sqsubset y \text{ for some }y \in a \}$
          is clopen for every clopen $a \subseteq X$.
  \end{enumerate}
  These constitute the category $\cat{SIS}$, whose morphisms are continuous
  morphisms that are bounded with respect to both relations.
\end{definition}

\begin{theorem}
  We have $\DFrm \cong \cat{SIS}$.
\end{theorem}
\begin{proof}
  Let $(X, \preceq, \sqsubset, P)$ be a descriptive $\sto$-frame and write
  $\tau_P$ for the topology on $X$ generated by the subbase $P$.
  Since $(X, \preceq, P)$ is a descriptive intuitionistic Kripke frame
  we know that $(X, \preceq, \tau_P)$ is an Esakia space.
  Furthermore, \ref{it:sis-1} follows from the fact that $(X, \preceq, \sqsubset)$
  is a $\sto$-frame and \ref{it:sis-3} follows from $\sqsubset$-refinedness.
  For \ref{it:sis-2}, let $a$ be a clopen subset of $X$.
  Then $a = \bigcup_{i = 1}^n b_i \cap -c_i$,
  where $b_i, c_i \in P$,
  and we have
  \begin{align*}
    {\downarrow}_{\sqsubset}a
      &= \bigcup {\downarrow}_{\sqsubset}(b_i \cap -c_i) \\
      &= \bigcup -(b_i \undsto c_i)
  \end{align*}
  which is in $-P$.
  
  Conversely, for a strict implication space $(X, \preceq, \sqsubset, \tau)$
  let $P$ be the collection of clopen upsets of $(X, \tau)$.
  Then $(X, \preceq, P)$ is a descriptive intuitionistic Kripke frame,
  $(X, \preceq, \sqsubset)$ is a $\sto$-frame because of \ref{it:sis-1},
  $\sqsubset$-refinedness follows from \ref{it:sis-3} and an argument that
  resembles the proof of Proposition, and closure of $P$ under $\undsto$
  follows from the fact that $a \undsto b = X \setminus {\downarrow}_{\sqsubset}(a \cap -b)$
  is in $P$ as a consequence of \ref{it:sis-2}.
  
  It is obvious that these two transformations define a bijection on objects.
  The isomorphism on morphisms is trivial.
\end{proof}

\subsection{Proof of Lemma \ref{lem:key}}\label{subsec:proof-lem-key}

  For the proof of Lem.~\ref{lem:key} we make use of the following
  auxiliary lemma.

\begin{lemma}\label{lem:sto}
  Let $(A, \sto)$ be an \HL-algebra.
  Then for all $a, b, c \in A$ we have
  \begin{enumerate}
    \item \label{it:lem:sto-1}
          $a \sto b \leq (a \wedge c) \sto b$ (antitone in first argument)
    \item \label{it:lem:sto-3}
          $a \sto b \leq (a \wedge c) \sto (b \wedge c)$
  \end{enumerate}
\end{lemma}
\begin{proof}
  It follows from \ref{it:C2} and absorption that
  $$
    (a \sto b) \wedge ((a \wedge c) \sto b)
      = (a \vee (a \wedge c)) \sto b = a \sto b,
  $$
  and this entails the first item.
  The second item follows from \cite[Lem.~4.1(a)]{LitVis18}.
\end{proof}

\begin{proof}[Proof of Lemma \ref{lem:key}]
  We prove two inclusions.
  Suppose $\mf{p} \in \widetilde{a \sto b}$, so that $a \sto b \in \mf{p}$.
  If $\mf{p} \sqsubset \mf{q}$ and $\mf{q} \in \tilde{a}$, then $a \in \mf{q}$.
  By definition of $\sqsubset$ this implies $b \in \mf{q}$, hence $\mf{q} \in \tilde{b}$.
  Therefore $\mf{p} \in \tilde{a} \undsto \tilde{b}$.
  
  For the converse we need to work harder.
  Suppose $\mf{p} \notin \widetilde{a \sto b}$.
  We aim to find a prime filter $\mf{q} \in \fun{pf}A$ such that:
  \begin{enumerate}[\quad (1)]
    \item \label{it:pf-q-1}
          $\mf{p} \sqsubset \mf{q}$;
    \item \label{it:pf-q-2}
          $a \in \mf{q}$;
    \item \label{it:pf-q-3}
          $b \notin \mf{q}$.
  \end{enumerate}
  By definition of $\undsto$ this entails
  $\mf{p} \notin \tilde{a} \undsto \tilde{b}$.
  
  We aim to construct the desired prime filter $\mf{q}$ using the prime
  filter lemma. To this end,
  let $F = \{ c \in A \mid a \preceq c \}$ %
  and $I = \{ d \in A \mid d \sto b \in \mf{p} \}$.
  Trivially $F$ is a filter, and it follows from \ref{it:C2} that $I$ is an ideal.
  In particular, we have $a \in F$ and $b \in I$.
  Moreover, $F \cap I = \emptyset$. To see this, suppose $x \in F \cap I$.
  Then $a \preceq x$ and $x \sto b \in p$, so by Lem.~\ref{lem:sto}\eqref{it:lem:sto-1}
  $x \sto b \preceq a \sto b$ and since $p$ is a filter $a \sto b \in \mf{p}$.
  A contradiction.
  
  Thus we can invoke the prime filter lemma to obtain a prime filter $\mf{q}$
  containing $F$ and disjoint from $I$.
  But more is true: if we have a look at the proof of the prime filter lemma
  (see e.g.~\cite[Lem.~1.4]{Mor05})
  we see that $\mf{q}$ is a maximal element
  with the property that it is a filter containing $F$ disjoint from $I$.
  That is, it is a maximal element in the set
  $$
    \ms{P} = \{ \mf{s} \in \fun{fil}A \mid F \subseteq \mf{s} \text{ and }
                                            \mf{s} \cap I = \emptyset \},
  $$
  which is ordered by inclusion.
  
  Items \eqref{it:pf-q-2} and \eqref{it:pf-q-3} are already satisfied
  We will now prove that $\mf{q}$ is such that $\mf{p} \sqsubset \mf{q}$.
  Suppose towards a contradiction that $\mf{p} \not\sqsubset \mf{q}$.
  Then we can find $c, d \in A$ such that $c \sto d \in \mf{p}$ and $c \in \mf{q}$ and
  $d \notin \mf{q}$. Let $\mf{q}_d$ be the filter generated by $\mf{q} \cup \{ d \}$.
  It is easy to check that
  $$
    \mf{q}_d = \{ c \wedge e \mid c \in \mf{q}, d \preceq e \}.
  $$
  Since $\mf{q}_d$ properly contains $\mf{q}$ and $\mf{q}$ is maximal in $\ms{P}$,
  we must have $\mf{q}_d \notin \ms{P}$.
  This implies $\mf{q}_d \cap I \neq \emptyset$, and therefore we can
  find $e \in \mf{q}$ such that $e \wedge d \in I$.
  By definition of $I$ this means $(e \wedge d) \sto b \in \mf{p}$.
  By assumption $c \sto d \in \mf{p}$ and so it follows from Lem.~\ref{lem:sto}\eqref{it:lem:sto-3}
  that $(e \wedge c) \sto (e \wedge d) \in \mf{p}$.
  By \ref{it:C3}, 
  $$
    \big((e \wedge c) \sto (e \wedge d)\big)
      \wedge \big((e \wedge d) \sto b \big)
      \preceq (e \wedge c) \sto b,
  $$
  and since $\mf{p}$ is a filter this implies $(e \wedge c) \sto b \in \mf{p}$.
  By definition of $I$ we now have $e \wedge c \in I$.
  But we chose both $e$ and $c$ from $\mf{q}$, and since $\mf{q}$ is a filter
  this implies $e \wedge c \in \mf{q}$.
  This is a contradiction with the fact that $\mf{q}$ and $I$ are disjoint,
  hence the assumption that $\mf{p} \not\sqsubset \mf{q}$ must be false.
  We conclude that $\mf{p} \sqsubset \mf{q}$.
  This completes the proof of the lemma.
\end{proof}

\begin{proof}[Proof of Proposition \ref{prop:iA-to-gfrm}]
  We know from intuitionistic logic that $(\fun{pf}A, \subseteq, \tilde{A})$
  is a descriptive intuitionistic Kripke frame.
  To show that $(\fun{pf}A, \subseteq, \sqsubset)$ is a $\sto$-frame,
  suppose $\mf{p}' \subseteq \mf{p} \sqsubset \mf{q}$.
  Then $a \sto b \in \mf{p}'$ implies $a \sto b \in \mf{p}$, so
  whenever $a \in \mf{q}$ we also have $b \in \mf{q}$.
  So $\mf{p}' \sqsubset \mf{q}$.
  Lemma \ref{lem:key} entails that $\tilde{A}$ is closed under
  $\undsto$ because
  $$
    \tilde{a} \undsto \tilde{b} = \widetilde{a \sto b} \in \tilde{A}.
  $$
  Finally, it follows from the construction of $\sqsubset$
  that $\alg{A}_*$ is $\sqsubset$-refined, and since it is based on a
  descriptive intuitionistic Kripke frame, $\alg{A}_*$ is descriptive.
\end{proof}

\begin{proof}[Proof of Lemma \ref{lem:nat-iso-iA}]
  This follows from the fact that the maps involved piggy-back on
  those from the duality between Heyting algebras and descriptive intuitionistic
  Kripke frames, which are natural.
\end{proof}

\begin{proof}[Proof of Theorem \ref{thm:duality}]
  It suffices to prove that we have natural isomorphisms
  $\fun{id}_{\HLAs} \cong ((\cdot)_*)^*$ and
  $\fun{id}_{\DFrm} \cong ((\cdot)^*)_*$.
  The former was proved in Lem.~\ref{lem:nat-iso-iA},
  so we focus on the latter.
  
  Let $\mf{G} = (X, \preceq, \sqsubset, P)$ be a descriptive frame.
  We already know that $\hat{(\cdot)} : \mf{G} \to (\mf{G}^*)_*$
  given by $\hat{x} = \{ a \in P \mid x \in a \}$ is a natural isomorphism
  between the underlying descriptive intuitionistic Kripke frames \cite[\S~8.4]{ChaZak97}.
  So it suffices to prove that $\hat{(\cdot)}$ is a $\sto$-frame morphism.
  Since $\hat{(\cdot)}$ is a bijection, it is enough to prove that
  $x \sqsubset y$ if and only if $\hat{x} \sqsubset \hat{y}$.
  But this follows immediately from the definition of our functors:
  \begin{align*}
    x \sqsubset y
      &\jff \forall a, b \in P \; [
            x \in a \undsto b \text{ and } y \in a \text{ imply } y \in b]  \\
      &\jff \forall a, b \in P \; [
            a \undsto b \in \hat{x} \text{ and } a \in \hat{y} \text{ imply } b \in \hat{y}] \\
      &\jff \hat{x} \sqsubset \hat{y}.
  \end{align*}
  This proves the theorem.
\end{proof}

\subsection{Proof of Lemma \ref{lem:translation-pres}}\label{subsec:translation-pres-proof}
  
  \noindent
  The proof of Lem.~\ref{lem:translation-pres} follows from the following
  two lemmas. Before each of these we introduce the relevant notion of the
  translation of a valuation.

  Let $\mf{F} = (X, \Ri, \Rm , P)$ be a general $\sf{S4K}$-frame.
  If $V : \Prop \to P$ is a valuation for it, then we define the valuation
  $\hrho V : \Prop \to \hrho P$ for $\hrho\mf{F} = (\eqc{X}, \eqc{\Ri}, \eqc{\Rm ^*}, \hrho P)$
  by
  $$
    \hrho V(p) := \eqc{\ibx V(p)} = \ibx \big(\textstyle\bigcup \eqc{\ibx V(p)}\big).
  $$
  
\begin{lemma}\label{lem:translation-pres-back}
  Let $\mf{F} = (X, \Ri, \Rm , P)$ be a general $\sf{S4K}$-frame and
  $V$ a valuation for $\mf{F}$. Then we have
  $$
    (\mf{F}, V), x \Vdash t(\phi) \jff (\hrho\mf{F}, \hrho V), \eqc{x} \Vdash \phi
  $$
  for all $x \in X$ and $\phi \in \lan{L}_{\sto}$.
\end{lemma}
\begin{proof}
  By induction on the structure of $\phi$.
  The cases $\top$ and $\bot$ are obvious.
  
  \medskip\noindent
  \fbox{$\phi = p \in \Prop$}
  We have
  \begin{align*}
    (\mf{F}, V), x \Vdash t(p) = \Boxi p
      &\jff \Ri[x] \subseteq V(p) \\
      &\jff x \in \ibx V(p) \\
      &\jff \eqc{x} \in \eqc{\ibx V(p)} = \hrho V(p) \\
      &\jff (\hrho\mf{F}, \hrho V), \eqc{x} \Vdash p
  \end{align*}
  The right-to-left direction from the third ``iff'' follows from the fact that
  $x \sim x'$ implies $x' \in \ibx V(p)$.
  
  \medskip\noindent
  \fbox{$\phi = \phi_1 \wedge \phi_2$}
  We have
  \begin{align*}
    (\mf{F}, V), x &\Vdash t(\phi_1 \wedge \phi_2) = \Boxi (t(\phi_1) \wedge t(\phi_2)) \\
      &\jff x\Ri y \To y \Vdash t(\phi_1) \text{ and } y \Vdash t(\phi_2) \\
      &\jff \eqc{x}\eqc{\Ri}\eqc{y} \To \eqc{y} \Vdash \phi_1 \text{ and } \eqc{y} \Vdash \phi_2 \\
      &\jff \eqc{x}\eqc{\Ri}\eqc{y} \To \eqc{y} \Vdash \phi_1 \wedge \phi_2 \\
      &\jff (\hrho\mf{F}, \hrho V), \eqc{x} \Vdash \phi_1 \wedge \phi_2
  \end{align*}
  The second ``iff'' follows from the induction hypothesis and the definition of $\eqc{\cdot}$.
  
  \medskip\noindent
  \fbox{$\phi = \phi_1 \vee \phi_2$}
  Similar to the previous case.
  
  \medskip\noindent
  \fbox{$\phi = \phi_1 \to \phi_2$}
  Compute
  \begin{align*}
    (\mf{F}, V), x &\Vdash t(\phi_1 \to \phi_2) = \Boxi (t(\phi_1) \to t(\phi_2)) \\
      &\jff x\Ri y \text{ and } y \Vdash t(\phi_1) \text{ imply } y \Vdash t(\phi_2) \\
      &\jff \eqc{x}\eqc{\Ri}\eqc{y} \text{ and } \eqc{y} \Vdash \phi_1 \text{ imply } \eqc{y} \Vdash \phi_2 \\
      &\jff (\hrho\mf{F}, \hrho V), \eqc{x} \Vdash \phi_1 \to \phi_2
  \end{align*}

  \medskip\noindent
  \fbox{$\phi = \phi_1 \sto \phi_2$}
  First assume $(\mf{F}, V), x \Vdash t(\phi_1 \sto \phi_2)$. Compute
  \begin{align*}
    (\mf{F}, \, & V), x \Vdash t(\phi_1 \sto \phi_2) = \Boxi \Boxm (t(\phi_1) \to t(\phi_2)) \\
      &\jff x\Ri y \text{ implies } y \Vdash \Boxm (t(\phi_1) \to t(\phi_2)) \\
      &\jff x(\Ri \circ \Rm )z \text{ implies } z \Vdash t(\phi_1) \to t(\phi_2) \\
      &\jff x(\Ri \circ \Rm )z \text{ and } z \Vdash t(\phi_1) \text{ imply } z \Vdash t(\phi_2) \\
      &\jff x\Rm ^*z \text{ and } z \Vdash t(\phi_1) \text{ imply } z \Vdash t(\phi_2)
  \end{align*}
  Now suppose $\eqc{x}\eqc{\Rm ^*}\eqc{z}$ and $\eqc{z} \Vdash \phi_1$.
  Then there exists $z' \in X$ such that $z \sim z'$ and $x\Rm ^*z'$.
  Since $\eqc{z} \Vdash \phi_1$ by the induction hypothesis we have $z' \Vdash t(\phi_1)$,
  so by the derivation above $z' \Vdash t(\phi_2)$ and therefore
  $\eqc{z} = \eqc{z'} \Vdash \phi_2$. So $\eqc{x} \Vdash \phi_1 \sto \phi_2$.
  
  Conversely, suppose $\eqc{x} \Vdash \phi_1 \sto \phi_2$.
  Then $x\Rm ^*z$ implies $\eqc{x}\eqc{\Rm ^*}\eqc{z}$ and the desired result follows
  from the induction hypothesis.
\end{proof}

  Let $\mf{F} = (X, \Ri, \Rm , P)$ be a general $\sf{S4K}$-frame.
  Suppose $W : \Prop \to \hrho P$ is a valuation for $\hrho\mf{F}$.
  Then since $W(p) \in \hrho P$ it must be of the form $\ibx \eqc{a}$ for some
  $a$ such that $\bigcup\eqc{a} \in P$. We pick such $a$ and call it $W'(p)$.

\begin{lemma}\label{lem:translation-pres-forw}
  Let $\mf{F} = (X, \Ri, \Rm , P)$ be a general $\sf{S4K}$-frame
  and $W$ a valuation for $\hrho\mf{F}$.Then we have
  $$
    (\mf{F}, W'), x \Vdash t(\phi) \jff (\hrho\mf{F}, W), \eqc{x} \Vdash \phi
  $$
  for all $x \in X$ and $\phi \in \mc{L}_{\sto}$.
\end{lemma}
\begin{proof}
  By induction on the structure of $\phi$.
  The cases $\top$ and $\bot$ are obvious.
  
  \medskip\noindent
  \fbox{$\phi = p \in \Prop$}
  If $(\mf{F}, W'), x \Vdash t(p) = \Boxi p$ then $x\Ri y$ implies $y \in W'(p)$.
  By definition of $\eqc{\cdot}$ we have $x\Ri y$ iff $\eqc{x}\eqc{\Ri}\eqc{y}$
  and by definition of $W'(p)$ we have $y \in W'(p)$ iff $\eqc{y} \in \eqc{W'(p)}$.
  Therefore $\eqc{x} \in \ibx\eqc{W'(p)} = W(p)$, so that
  $(\hrho\mf{F}, W), \eqc{x} \Vdash p$.
  
  Conversely, suppose $(\hrho\mf{F}, W), \eqc{x} \Vdash p$.
  Then $\eqc{x} \in \ibx\eqc{W'(p)}$,
  so $\eqc{x}\eqc{\Ri}\eqc{y}$ implies $\eqc{y} \in \eqc{W'(p)}$.
  Again as a consequence of the definitions this gives
  $x\Ri y$ implies $y \Vdash p$, so that $(\mf{F}, W'), x \Vdash \Boxi p = t(p)$.

  \medskip\noindent
  \fbox{All other cases}
  are the same as in Lem.~\ref{lem:translation-pres-back}.
\end{proof}

\subsection{Proof of Theorem \ref{thm:fmp-total}}\label{app:proof-thm-fmp-total}

  The proof of Thm.~\ref{thm:fmp-total} follows from the following
  two theorems, which consider subframe logics and $\Rm$-cofinal subframe logics,
  respectively.

\begin{theorem}\label{thm:fmp-easy}
  Suppose $\Theta$ is a canonical extension of $\lna{S4} \otimes \lna{K4}$
  containing \ref{ax:HL} that is closed under forming subframes. Then:
  \begin{enumerate}
    \item \label{it:thm:bm-1}
          $\Theta$ has the finite model property.
    \item \label{it:thm:bm-2}
          If moreover $\Theta$ contains the classical strength axiom
\begin{enumerate} 
 \axiom{S_c} \label{ax:Sc}
    $\Boxi p \to \Boxm p$.
 \end{enumerate}   
          then for any subframe logic $\Gamma \subseteq \lan{L}_{\sf{m}}$,
          the logic $\Theta \oplus \Gamma$ has the finite model property.
  \end{enumerate}
\end{theorem}
\begin{proof}[Proof of Theorem \ref{thm:fmp-easy}]
  Let $\mf{F} = (X, \Ri, \Rm, P)$ be a descriptive frame for $\Theta$ refuting a $\lanBM$-formula
  $\phi$ under valuation $V$.
  Then we will construct a finite subframe of $\mf{F}$ that validates $\Theta$
  and refutes $\phi$.
  In order to prove this, we an adaptation of the proofs of Theorems 17 and 21
  in \cite{WolZak98}. 
  
  We begin by constructing an inductive sequence $\{ X_i \}_{i \in \omega}$
  of subsets of $X$.
  
  \medskip\noindent
  {\it\underline{Base step.}}
  By \cite[Lem.~14]{WolZak98} we can pick a $\Ri$-maximal state $x_0$ where
  $\phi$ is refuted. Let $X_0 = \{ x_0 \}$.
  
  \medskip\noindent
  {\it\underline{Odd inductive step.}}
  Suppose $n$ is even. For each $x' \in X_n$ pick an $\Ri$-maximal witness
  for each $\phi$-equivalence class above $x'$ in the $\Ri$-order, and denote the
  set of such witnesses by $W_{x'}$. 
  (We can find such witnesses by \cite[Lem.~14]{WolZak98}.)
  Let
  $$
    X_{n+1} = X_n \cup \bigcup_{x' \in X_n} W_{x'}.
  $$
  The sets $W_{x'}$ are finite because there are only finitely many
  $\phi$-equivalence classes.
  Note also that in the presence  of the strength axiom \ref{ax:Sa},
  $\Ri$-maximal successors are also $\Rm$-maximal ones. 
  
  \medskip\noindent
  {\it\underline{Even inductive step.}}
  Suppose $n$ is odd.
  For each $x' \in X_{n+1}$, using \cite[Lem.~14]{WolZak98} 
  pick an $\Rm$-maximal witness for each
  $\phi$-equivalence class, and denote the set of such witnesses by $W_{x'}$.
  Again, define
  $
    X_{n+1} = X_n \cup \bigcup_{x' \in X_n} W_{x'}.
  $

  \medskip\noindent
  {\it\underline{The set $X_{\omega}$.}}
  Define
  $$
    X_{\omega} := \bigcup_{n \in \omega} X_n.
  $$%
  Before turning this into a model, we shall argue that it is a finite
  subset of $X$. Our more general setting compared to \cite{WolZak98} makes
  this more complicated, but not impossible.
  
  For each $x \in X_{\omega}$, denote by $n(x)$ the smallest integer
  $n \in \omega$ such that $x \in X_n$.
  Define relations $\Ri'$ and $\Rm'$ on $X_{\omega}$ by
  \begin{align*}
    x \Ri' y &\quad\text{ if }\quad x \Rm y \text{ and } n(y) \text{ is odd and } y \in X_{n(y)-1} \\
    x \Rm' y &\quad\text{ if }\quad x \Rm y \text{ and } n(y) \text{ is even and } y \in X_{n(y)-1}
  \end{align*}
  Then $x \Ri' y$ implies that $y$ is $\Ri$-maximal relative to $\sim_{\phi}$,
  and similar for $\Rm'$.
  Clearly the structure $(X_{\omega}, \Ri' \cup \Rm')$, viewed as a graph,
  is connected.
  Therefore we can invoke K\"{o}nig's Lemma to obtain an infinite sequence.
  
  If this sequence contains an infinite number of $\Rm'$ transitions then
  there exists an infinite subsequence of the form
  $$
    x_0 \Rm' x_1 \Rm' x_2 \Rm' x_3 \cdots
  $$
  By construction each of the $x_i$ is $\Rm$-maximal.
  But then transitivity of $\Rm$ implies that each of the $x_i$
  belong to a different $\sim_{\phi}$-equivalence class.
  A contradiction, since there are only finitely many such classes.
  
  If the sequence obtained from K\"{o}nig's Lemma has a finite number of
  $\Rm'$-transitions, then there must be an infinite subsequence of the form
  $
    x_0 \Ri' x_1 \Ri' x_2 \Ri' x_3 \cdots
  $
  and a similar argument as above yields a contradiction.
  Thus no infinite sequence can exist, and therefore $X_{\omega}$ must be finite.
  
  We could not simply invoke K\"{o}nig's lemma to $X_{\omega}$ ordered by
  (restrictions of) $\Ri$ and $\Rm$ (like in \cite{WolZak98}), because this
  could potentially yield an infinite chain whose states are not all
  $\Ri$-maximal or $\Rm$-maximal.

  \medskip\noindent
  {\it\underline{Finite submodel of $\mf{F}$.}}
  Define the frame $\mf{F}'$ to be the subframe of $\kappa\mf{F}$
  generated by $X_{\omega}$,
  and $\mf{M}'$ as the model $(\mf{F}', V')$, where $V'(p) = V(p) \cap (X_0 \cup X_{\omega})$.
  Then by Lem.~\ref{lem:subsub-truth} we have
  $$
    \mf{M}, y \Vdash \phi \jff \mf{M}', y \Vdash \psi
  $$
  for all $y \in X_{\omega}$ and $\psi \in \Subf(\phi)$.
  Therefore $\mf{F}'$ refutes $\phi$.
  Since $\Theta$ is a canonical subframe logic and,
  $\mf{F}'$ is a subframe of $\kappa\mf{F}$, we also have $\mf{F}' \Vdash \Theta$.
  This proves item \ref{it:thm:bm-1}.
  
  \medskip\noindent
  {\it Proof of item \ref{it:thm:bm-2}.}
  Assume the above construction started with a descriptive frame $\mf{G}$
  for $\Theta \oplus \Gamma$ that refutes $\phi$.
  We have already seen that the resulting frame $\mf{F}'$ refutes $\phi$
  and validates $\Theta$, so it remains to show that $\mf{F}' \Vdash \Gamma$.
  
  The assumption of strength, together with the construction of $X_{\omega}$,
  implies that each state $y \in X_{\omega}$ is $\Rm$-maximal.
  Therefore, since $\Rm$ is transitive, it follows from \cite[Lem.~15]{WolZak98}
  that $\mf{G}' \Vdash \Gamma$.
\end{proof}

\begin{theorem}\label{thm:fmp-cofinal}
  Suppose $\Theta$ is a canonical extension of $\lna{S4} \otimes \lna{K4}$
  containing \ref{ax:HL} that is closed under forming $\Rm$-cofinal subframes. Then:
  \begin{enumerate}
    \item \label{it:thm:bm-1}
          $\Theta$ has the finite model property.
    \item \label{it:thm:bm-2}
          If $\Theta$ contains the strength axiom \ref{ax:Sa},
          then for any $\Rm$-subframe logic $\Gamma \subseteq \lan{L}_{\sf{m}}$,
          the logic $\Theta \oplus \Gamma$ has the finite model property.
  \end{enumerate}
\end{theorem}
\begin{proof}[Proof of Theorem \ref{thm:fmp-cofinal}]
  Let $\mf{F} = (X, \Ri, \Rm, P)$ be a descriptive frame for $\Theta$
  that refutes $\phi$.
  Let $V$ be a valuations such that $(\mf{F}, V) \not\Vdash \phi$.
  In order to prove the theorem, we modify the proof of
  Thm.~\ref{thm:fmp-easy} as follows:
  First, we modify the construction to obtain an $\Rm$-cofinal subframe.
  Second, we quotient out this new subframe to make it finite.

  Let $\{ C_j \mid j \in J \}$ be the set of all $\Rm$-final $\Rm$-clusters.
  Since $\mf{F}$ is descriptive, every state $x$ has an $\Rm$
  successor in an $\Rm$-final $\Rm$-cluster \cite[Thm.~10.36]{ChaZak97}.
  Therefore, to achieve $\Rm$-cofinality, it suffices to add to $X_{\omega}$
  a state from each $C_j$.
  In order to still be able to use Lem.~\ref{lem:subsub-truth} we add a
  (finite) set $F_j$ of states satisfying the precondition from
  Lem.~\ref{lem:subsub-truth} for each $j \in J$.
  
  \medskip\noindent
  {\it\underline{Constructing $F_j$.}}
  Let $C_j$ be a $\Rm$-final $\Rm$-cluster.
  Let $F_j'$ be a minimal subset of $C_j$ such that for each $x \in C_j$
  there is an $\Rm$-maximal $y \in F_{j,0}$ such that $x \sim_{\phi} y$.
  (We can find such $\Rm$-maximal states using \cite[Lem.~14]{WolZak98}.)
  Since $C_j$ is a cluster and there are only finitely many $\sim_{\phi}$-equivalence
  classes the set $F_j'$ is finite.
  
  Now suppose $F_{j,k}$ has been defined. We give $F_{j,k+1}$.
  For each $y \in F_{j,k}$ let $Y_y$ be a minimal set of maximal $\Ri$-states
  $\Ri$-above $y$, such that for ever $z$ with $y \Ri z$ there exists $z' \in Y_y$
  such that $z \sim_{\phi} z'$.
  Define
  $$
    F_{j,k+1} = \bigcup_{y \in F_{j,k}} Y_y.
  $$
  We claim that this process is finite. If $z \in F_{j,k}$ is introduced in
  an earlier step, then $Y_z$ (used in the construction of $F_{j,k+1}$) is
  empty by minimality. If $z$ was introduced in the construction of
  $F_{j,k}$ as an element of some $Y_y$, where $y \in F_{j,k-1}$,
  then by construction and $\Ri$-maximality of the states in $F_{j,k-1}$
  there are at most $c - k$ different $\sim_{\phi}$-equivalence classes
  that $z$ can see. (Recall that $c$ denotes the number of $\sim_{\phi}$-equivalence
  classes.) Therefore $|Y_z| \leq c - k$. This proves that the recursion terminates
  after $c$ steps.
  
  Therefore, the set
  $$
    F_j = \bigcup_{0 \leq k \leq c} F_{j,k}
  $$
  is finite. In fact, each of the $F_j$ is bounded by $c^{c+1}$.
  Moreover, we claim that it satisfies the precondition from Lem.~\ref{lem:subsub-truth}.
  Clearly, if $y \in F_j$ and $y \Ri z$, then by construction there exists
  $z' \in F_j$ such that $z \sim_{\phi} z'$.
  If $y \in F_j$ and $y \Rm z$, then by \eqref{eq:p-sto} we have $z \in C_j$.
  By construction there exists $z' \in C_j$ such that $z \sim_{\phi} z'$
  and by transitivity of $\Rm$ this implies $y \Rm z'$.
  
  \medskip\noindent
  {\it\underline{An $\Rm$-cofinal subframe.}}
  Let $X_{\omega}$ be constructed as in the proof of Thm.~\ref{thm:fmp-easy}.
  Then by construction the subframe $\mf{F}^{\dagger}$ of $\mf{F}$ generated by
  $$
    \ov{X}_{\omega} = X_{\omega} \cup \bigcup \{ F_j \mid j \in J \}
  $$
  is $\Rm$-cofinal. Let $V^{\dagger}$ be the induced valuation and set
  $\mf{M}^{\dagger} = (\mf{F}^{\dagger}, V^{\dagger})$.
  Then since $\Theta$ is closed under $\Rm$-cofinal subframes we have
  $\mf{F} \Vdash \Theta$, and as a consequence of Lem.~\ref{lem:subsub-truth}
  we have $\mf{M}^{\dagger} \not\Vdash \phi$, so that $\mf{F}^{\dagger}\not\Vdash \phi$.
  
  \medskip\noindent
  {\it\underline{Making $\mf{M}^{\dagger}$ finite.}}
  We use a trick similar to \cite[Thm.~21]{WolZak98}.
  For each $j \in J$, let $F_j^+ = (C_j \cap X_{\omega}) \cup F_j$.
  Since the $F_j$ and $X_{\omega}$ are uniformly bounded by $c^{c+1} + |X_{\omega}|$,
  so are the sets $F_j^+$ are uniformly bounded.
  Hence there are only finitely many non-isomorphic submodels of $\mf{M}^{\dagger}$
  generated by $F_j^+$.
  Identifying isomorphic such submodels yields a finite quotient
  of $\mf{M}^{\dagger}$ that witnesses the claim.
  
  \medskip\noindent
  {\it\underline{Item 2.}}
  The second item is proved in a similar way as the second item of
  Thm.~\ref{thm:fmp-easy}.
\end{proof}

\end{document}